\documentclass[12pt]{amsart}
\usepackage{a4wide}
\usepackage{eucal}
\usepackage{umlaut}
\usepackage{latexsym}
\usepackage{amssymb}
\usepackage{verbatim}
\usepackage{graphicx}
\usepackage[active]{srcltx}
\usepackage[all]{xy}
\usepackage[pdfauthor   = {Mohamed Barakat \& Daniel Robertz},
            pdftitle    = {conley: Computing connection matrices in Maple},
            pdfsubject  = {18E30; 37B30},
            pdfkeywords = {Connection matrix; Conley index; Morse decomposition},
            bookmarks=true,
            bookmarksopen=true,
            colorlinks=true,
            pagebackref=true,
            hyperindex=true,
            linkcolor=blue,
            pagecolor=blue,
            citecolor=blue,
            urlcolor=blue,
            hypertex=true
            ]{hyperref}
\usepackage{maple2e}
\LeftMapleSkip=0em

\usepackage{mathdots}

\input amathmoe.sty



\DeclareMathOperator{\Aut}{Aut}

\renewcommand\phi{\varphi}
\newcommand{\img}{\mathrm{im}}
\newcommand{\Inv}{\mathrm{Inv}}

\newcommand\homalg{\mathtt{homalg}}
\newcommand\conley{\mathtt{conley}}
\newcommand\Maple{\mathsf{Maple}}
\newcommand\NN{\mathcal{N}}

\newcommand{\id}{\mathrm{id}}
\newcommand{\Con}{\mathrm{Con}}


\bdoc

\author{Mohamed Barakat \& Daniel Robertz}
\address{Lehrstuhl B f\"ur Mathematik, RWTH-Aachen University, 52062 Germany}
\email{mohamed.barakat@rwth-aachen.de, daniel@momo.math.rwth-aachen.de}

\PUSH{conley_package.tex}%
\title{$\mathtt{conley}$\\ \medskip Computing connection matrices in $\Maple$}

\begin{abstract}
  In this work we announce the $\Maple$ package $\conley$ to compute connection and $C$-connection matrices. $\conley$ is based on our abstract homological algebra package $\homalg$. We emphasize that the notion of braids is irrelevant for the definition and for the computation of such matrices. We introduce the notion of triangles that suffices to state the definition of ($C$)-connection matrices. The notion of octahedra, which is equivalent to that of braids is also introduced.
\end{abstract}

\maketitle

\tableofcontents

\section{Introduction}

The algebraic theory of connection and $C$-connection matrices is concerned with condensing the information of a graded module octahedron in a matrix called a connection matrix or $C$-connection matrix.

In Section~\ref{Con} we recall the central notions of this algebraic theory, introduce some new notions and give a {\em braid-free} definition of ($C$)-connection matrices. For this we introduce the notion of a graded module triangular system or simply a graded module triangle. This is the only notion used in the definition of connection and $C$-connection matrices.

The role of braids is emphasized in formulating the necessary and sufficient conditions to the existence of $C$-connection matrices in Section~\ref{exist}. There we also introduce the equivalent notion of octahedra.

The notion of symmetric $C$-connection matrix is defined in
Section~\ref{symm}.

In Section~\ref{Morse} we briefly review the definition of the homology {\sc Conley} index of an isolated invariant set and the theory of {\sc Morse} decompositions of a dynamical system. This theory is the main source of applications of connection matrices.

The package $\conley$ is more or less a tiny application of our more elaborate homological algebra package $\homalg$ \cite{BR}. In Section~\ref{package} some technical comments on $\conley$ are made.

Finally, in Section~\ref{beispiele} we present some simple examples using $\conley$.

We are very much indebted to {\sc Stanislaus Maier-Paape}, who initiated the joint {\sc Conley} index seminar in Aachen, introduced us to the subject and explained us the fascinating
dynamical side of the theory. We would also like to thank him for fruitful discussions
which revealed a serious mistake in the first version of this paper. We also would like to
thank all the participants of the seminar, and especially {\sc Max Neunhoeffer} and
{\sc Felix Noeske} for discussions about the algebraic aspects of the theory.

Aspects of dynamical systems are not in the focus of this paper.
However, joint work with {\sc Maier-Paape} on more substantial dynamical applications and in particular on transition matrices \cite{Mis,MMW} is in preparation.
Transition matrices can also be computed using $\conley$. We don't touch upon this theory here.

Note that we apply morphisms from the {\em right} and hence we use the {\em row convention} for matrices. As one consequence we talk about {\em lower} triangular instead of upper triangular matrices.

\section{Connection Matrices}\label{Con}
\subsection{Posets}
A set $P$ together with a strict partial order $>$ (i.e.\ an irreflexive and transitive relation $> \, \subset P \times P$) is called a {\em poset}
and is denoted by $(P, >)$.

A subset $I \subset P$ is called an {\em interval}
in $(P, >)$
if for all $p, q \in I$ and $r \in P$ the following
implication holds:
\[
q > r > p \quad \Rightarrow \quad r \in I.
\]
The set of all intervals in $(P, >)$ is denoted by
${\mathcal I}(P, >)$.

An $n$-tuple $(I_1, \ldots, I_n)$ of intervals
in $(P, >)$ is called {\em adjacent} if these intervals are
mutually disjoint, $\bigcup_{i=1}^n I_i$ is an interval in
$(P, >)$ and for all $p \in I_j$, $q \in I_k$
 the following implication holds:
\[
j < k \quad \Rightarrow \quad p \not> q.
\]
The set of all adjacent $n$-tuples of intervals in $(P, >)$ is
denoted by ${\mathcal I}_n(P, >)$.

If $(I_1, \ldots, I_n)$ is an adjacent $n$-tuple of
intervals in $(P, >)$, then set
$I_1 I_2 \ldots I_n :=
\bigcup_{i=1}^n I_i$.

If $(I, J) \in {\mathcal I}_2(P, >)$ as well as
$(J, I) \in {\mathcal I}_2(P, >)$, then $I$ and $J$ are said to
be {\em noncomparable}.

\subsection{Triangles}
For notions like quasi-isomorphism and distinguished triangles we refer to \cite{GM}.

\begin{defn}[Graded module triangle]
Let $(P, >)$ be a poset. A {\em graded module triangular system} or simply a {\em graded module triangle} $G$ over $(P, >)$ consists of graded modules $G(I)$, $I \in {\mathcal I}(P, >)$, and homomorphisms of graded modules
\begin{eqnarray*}
i(I, IJ): G(I) & \longrightarrow & G(IJ),\\
p(IJ, J): G(IJ) & \longrightarrow & G(J),\\
\partial(J, I): G(J) & \longrightarrow & G(I)
\end{eqnarray*}
for all $(I, J) \in {\mathcal I}_2(P, >)$, where $i(I, IJ)$ and $p(IJ, J)$ are of degree $0$, $\partial(J, I)$ is of degree $-1$, and the following two conditions are satisfied:
\begin{enumerate}
\item[a)] For all $(I,J)\in{\mathcal I}_2(P, >)$ the triangle
\[
  \xymatrix@=1.5pc{
    & G(IJ)\ar[dr]^{p(IJ, J)} & \\
    G(I)\ar[ur]^{i(I, IJ)} && G(J)\ar[ll]^{\partial(J, I)}
  }
\]
is distinguished.
\item[b)] If $I, J \in {\mathcal I}(P, >)$ are noncomparable, then $i(I, IJ)p(JI, I) = \id_{|G(I)}$.
\end{enumerate}
\end{defn}

We think of a graded module as an element of $\mathrm{Kom}_0(\mathcal{A}) \subset\mathrm{Kom}(\mathcal{A})$.  $\mathrm{Kom}(\mathcal{A})$ is the category of complexes over an abelian category $\mathcal{A}$ and $\mathrm{Kom}_0(\mathcal{A})$ the complete subcategory of cyclic complexes, i.e.\  complexes with zero boundary maps. Cf.\  \cite[{III.2.3}]{GM}.

The one dimensional unraveling of condition a) is the following long exact sequence condition:
\begin{enumerate}
  \item [a')] For all $(I,J)\in{\mathcal I}_2(P, >)$ 
  \[
    \cdots \to G(I) \xrightarrow{i(I, IJ)} G(IJ)
    \xrightarrow{p(IJ, J)} G(J)
    \xrightarrow{\partial(J, I)} G(I) \to \cdots
\]
is a long exact sequence.
\end{enumerate}

\begin{defn}[Homomorphism of graded module triangles]
Let $(P, >)$ be a poset and let $G$ and $G'$ be graded module triangles over
$(P, >)$.
\begin{enumerate}
\item A {\em homomorphism of graded module triangles} $\theta: G \to G'$ consists of homomorphisms of graded modules $\theta(I): G(I) \to G'(I)$, $I\in\mathcal{I}(P,>)$, such that the following diagram commutes for all $(I, J) \in {\mathcal I}_2(P, >)$:
\begin{equation}\label{Isom}
  \xymatrix@=1.5pc{
    & G(IJ) \ar[dr]^{p(IJ, J)}\ar@{-}[d]|{\vphantom{\int^{x^2}}\theta(IJ)} & \\
    G(I)\ar[dd]_{\theta(I)}\ar[ur]^{i(I, IJ)} & \ar[d] &
    G(J)\ar[ll]^>(.4){\partial(J, I)}\ar[dd]^{\theta(J)} \\
    & G'(IJ)\ar[dr]^{p'(IJ, J)} & \\
    G'(I)\ar[ur]^{i'(I, IJ)} && G'(J)\ar[ll]^>(.4){\partial'(J, I)}
  }
\end{equation}
\item An {\em isomorphism of graded module triangles} $\theta: G \to G'$ is a homomorphism of graded module triangles where all $\theta(I)$ are isomorphisms of graded modules.
\end{enumerate}
\end{defn}

The two dimensional unraveling of the above diagram is:
\begin{equation}\label{isom}
\xymatrix@C=1.5pc@R=0.8pc{
\cdots \ar[r]
&
G(I) \ar[rr]^{i} \ar[dd]^{\theta(I)}
&
&
G(IJ) \ar[rr]^{p} \ar[dd]^{\theta(IJ)}
&
&
G(J) \ar[rr]^{\partial} \ar[dd]^{\theta(J)}
&
&
G(I) \ar[r] \ar[dd]^{\theta(I)}
&
\cdots
\\
&
&
&
&
&
&
&
&
\\
\cdots \ar[r]
&
G'(I) \ar[rr]^{i'}
&
&
G'(IJ) \ar[rr]^{p'}
&
&
G'(J) \ar[rr]^{\partial'}
&
&
G'(I) \ar[r]
&
\cdots
}
\end{equation}
\begin{defn}[Chain complex triangle]
Let $(P, >)$ be a poset. A {\em chain complex triangular system} or simply a {\em chain complex triangle} $C$ over $(P, >)$ consists of chain complexes $C(I)$, $I \in {\mathcal I}(P, >)$, and chain maps $i(I, IJ): C(I) \to C(IJ)$ and $p(IJ, J): C(IJ) \to C(J)$ for all $(I, J) \in {\mathcal I}_2(P, >)$, satisfying the following two conditions:
\begin{enumerate}
\item[a)] The complex $0 \to C(I) \xrightarrow{i(I, IJ)} C(IJ) \xrightarrow{p(IJ, J)} C(J) \to 0$ is quasi-isomorphic\footnote{More precisely, isomorphic in the localized category, in which one inverts quasi-isomorphisms, where quasi-isomorphisms are chain maps inducing isomorphism on homology \cite[{III.2}]{GM}.} to a distinguished triangle of complexes.
\item[b)] If $I, J \in {\mathcal I}(P, >)$ are noncomparable, then $i(I, IJ)p(JI, I) = \id_{|C(I)}$.
\end{enumerate}
\end{defn}

Short exact and short weakly exact sequences are examples of complexes quasi-isomorphic to a distinguished triangle of complexes. A sequence $0\to C'\xrightarrow{i}C\xrightarrow{p}C''\to 0$ of complexes is called short exact if $i$ is injective, $ip=0$, and $C/\img(i)\to C''$ induces isomorphism on homology. Using the cylinder-cone-translation construction \cite[{III.3}]{GM} one can construct distinguished triangles out of such sequences.

\begin{rmrk}[Chain complex generated triangles]
Let $(P, >)$ be a poset and $C$ a chain complex triangle over $(P, >)$.
The homology modules $H(I)$, $I \in {\mathcal I}(P, >)$, together with
the homomorphisms $\iota(I, IJ): H(I) \to H(IJ)$ and
$\pi(IJ, J): H(IJ) \to H(J)$ induced by the chain maps
$i(I, IJ)$ and $p(IJ, J)$ form a graded module triangle $HC$ over $(P, >)$.
\end{rmrk}

\subsection{$C$-Connection and connection matrices}

We fix a poset $(P, >)$. In what follows, we consider collections
$C=\{ C(p) \mid p \in P \}$ of graded modules, which are indexed by $P$,
and a homomorphism $\Delta: \bigoplus_{p \in P} C(p) \to \bigoplus_{p \in P} C(p)$.

For an interval $I$ in $(P, >)$
set $C(I) := \bigoplus_{p \in I} C(p)$ and denote by
$\Delta(I)$ the homomorphism
$\iota_I  \Delta  \pi_I$, where
$\iota_I: C(I) \to C(P)$ is the canonical injection and
$\pi_I: C(P) \to C(I)$ is the canonical projection.

If $p_1, p_2 \in P$, we refer to the restriction of $\Delta$ to $C(p_1)$ by $\Delta_{p_1}:C(p_1)\to C(P)$, and the composition $\Delta_{p_1}\pi_{p_2}$,
where $\pi_{p_2}$ is the projection $C(P) \to C_{p_2}$, is denoted by $\Delta_{p_1,p_2}:C(p_1)\to C(p_2)$. Then $\Delta$ can be visualized as a matrix with $\Delta_{p_1}$ as its $p_1$-th row and $\Delta_{p_1,p_2}$ as its entry at position $(p_1, p_2)$. The same applies to $\Delta(I)$ for $I \in {\mathcal I}(P, >)$, if $p_1, p_2 \in I$.

\begin{defn}[{\cite[Def.~1.3]{Fra88}}]
$\Delta$ being as above:
\begin{enumerate}
  \item $\Delta$ is said to be {\em lower triangular} if
$\Delta_{p_1,p_2} \neq 0$ implies $p_1 > p_2$ or $p_1 = p_2$.
  \item $\Delta$ is said to be {\em strictly lower triangular} if
$\Delta_{p_1,p_2} \neq 0$ implies $p_1 > p_2$.
  \item $\Delta$ is called a {\em boundary map} if each
$\Delta_{p_1,p_2}$ is a homomorphism of graded modules of degree $-1$ and
$\Delta \circ \Delta = 0$.
  \end{enumerate}
\end{defn}

\begin{prop}[{\cite[Prop.~3.3]{Fra89}}]
Let $C=\{ C(p) \mid p \in P \}$ be a collection of graded modules indexed by $P$ and let
$\Delta: \bigoplus_{p \in P} C(p) \to \bigoplus_{p \in P} C(p)$ be a
lower triangular boundary map. Then:
\begin{enumerate}
\item $C(I)$ and $\Delta(I)$ form a
chain complex $C^\Delta(I)$ for all $I \in {\mathcal I}(P, >)$.
\item For all $(I, J) \in {\mathcal I}_2(P, >)$,
the obvious injection and projection maps $i(I, IJ)$ and $p(IJ, J)$
are chain maps and
\[
0 \to C^\Delta(I) \xrightarrow{i(I, IJ)} C^\Delta(IJ) \xrightarrow{p(IJ, J)} C^\Delta(J) \to 0
\]
is a short exact sequence.
\end{enumerate}
\end{prop}

Using the above introduced notion the previous proposition simply asserts that $C^\Delta$ is a chain complex triangle.

The following definition of a connection matrix avoids braids (cf. {\cite[Def.~3.6]{Fra89}}).

\begin{defn}[($C$)-Connection matrix]\label{conn}
Let $G$ be a graded module triangle over $(P, >)$, $C=\{ C(p) \mid p \in P \}$ a collection of graded modules indexed by $P$, and $\Delta: \bigoplus_{p \in P} C(p) \to \bigoplus_{p \in P} C(p)$ a lower triangular boundary map. Let $C^\Delta$ be the chain complex triangle formed by the $C(I)$, $\Delta(I)$, $i(I, IJ)$, and $p(IJ, J)$.
\begin{enumerate}
  \item If $HC^\Delta$ and $G$ are isomorphic as graded module triangles, then
    $\Delta$ is called a {\em $C$-connection matrix} for $G$.
  \item If furthermore $C(p)\cong G(p)$, then $\Delta$ is called a {\em connection matrix} for $G$.
\end{enumerate}
\end{defn}

\section{Octahedra and Braids and the Existence of $C$-Connection Matrices}
\label{exist}

In this section we want to review the necessary and sufficient conditions a graded module triangle $G$ must satisfy so that a connection matrix (resp.\ $C$-connection matrix) for $G$ exists.

\subsection{Graded module octahedra}

\begin{defn}[Graded module octahedron]\label{octahedron}
Let $(P, >)$ be a poset. A {\em graded module octahedral system} or simply a {\em graded module octahedron} $G$ over $(P, >)$ is a graded module triangle over $(P,>)$ satisfying the extra condition:
\begin{enumerate}
  \item[c)]\label{octa} For all $(I, J, K) \in {\mathcal I}_3(P, >)$ the following
    octahedron commutes and is distinguished:
    
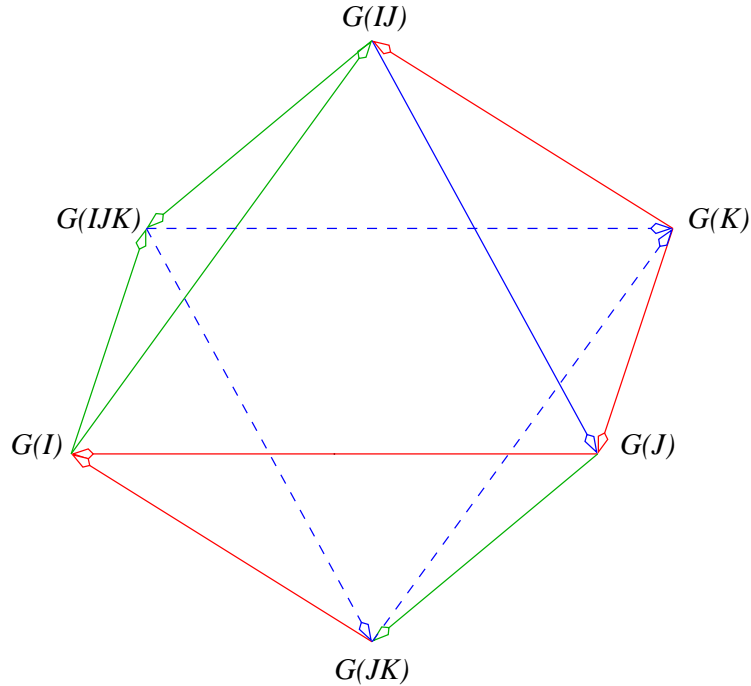
\begin{center}
      \includegraphics{octahedron.pstex}
    \end{center}
\end{enumerate}
\end{defn}

By ``commutes and is distinguished'' the following is meant: First notice that there are two types of triangles resp.\  squares appearing in the octahedron: cyclic and noncyclic ones. For the cyclic ones (4 triangles and one square) we require distinguishedness and for the others (4 triangles and two squares) commutativity. Actually the distinguishedness of the square, which is the horizontal one, follows from the rest.

The two dimensional unraveling of the octahedron condition c) is the following braid condition {\cite[Def.~1.1]{Fra88}}:
\begin{enumerate}
  \item[c')]\label{braid} For all $(I, J, K) \in {\mathcal I}_3(P, >)$ the following
braid diagram commutes:
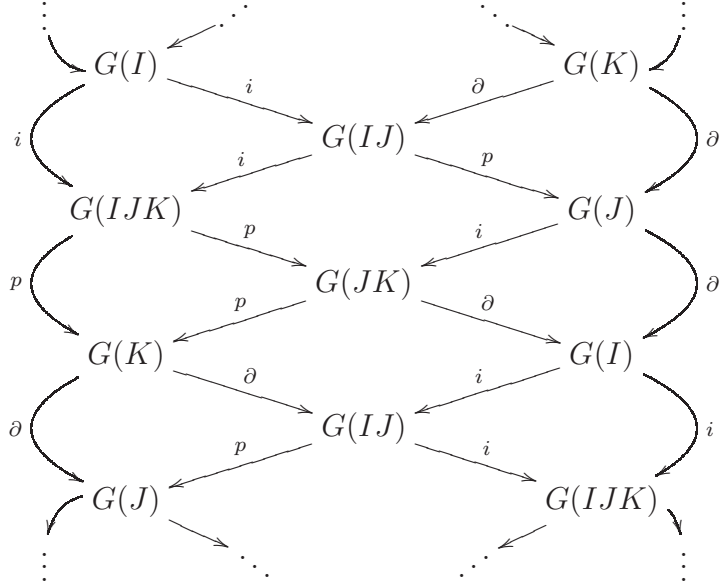
\[
\xymatrix@C=1.5pc@R=0.7pc{
\save[]-<6ex,0ex>*{\raisebox{0.5ex}{$\vdots$}}\ar@/_1pc/[d]\restore
&
\save[]-<-0.1ex,0ex>*{\raisebox{0.5ex}{$\iddots$} \; \; \; \;}\ar[dl]\restore
&
&
\save[]-<0.1ex,0ex>*{\; \; \; \; \raisebox{0.5ex}{$\ddots$}}\ar[dr]\restore
&
\save[]-<-6ex,0ex>*{\raisebox{0.5ex}{$\vdots$}}\ar@/^1pc/[d]\restore
\\
G(I) \ar[drr]^{i} \ar@/_3pc/[dd]_{i}
&
&
&
&
G(K) \ar[dll]_{\partial} \ar@/^3pc/[dd]^{\partial}
\\
&
&
G(IJ) \ar[dll]_{i} \ar[drr]^{p}
&
&
\\
G(IJK) \ar[drr]^{p} \ar@/_3pc/[dd]_{p}
&
&
&
&
G(J) \ar[dll]_{i} \ar@/^3pc/[dd]^{\partial}
\\
&
&
G(JK) \ar[dll]_{p} \ar[drr]^{\partial}
&
&
\\
G(K) \ar[drr]^{\partial} \ar@/_3pc/[dd]_{\partial}
&
&
&
&
G(I) \ar[dll]_{i} \ar@/^3pc/[dd]^{i}
\\
&
&
G(IJ) \ar[dll]_{p} \ar[drr]^{i}
&
&
\\
G(J)
&
&
&
&
G(IJK) \\
\save[]-<6ex,0ex>*{\vdots}\ar@{<-}@/^1pc/[u]\restore
&
\save[]-<-0.1ex,0ex>*{\; \ddots}\ar@{<-}[ul]\restore
&
&
\save[]-<0.1ex,0ex>*{\iddots \;}\ar@{<-}[ur]\restore
&
\save[]-<-6ex,0ex>*{\vdots}\ar@{<-}@/_1pc/[u]\restore
}
\]
\end{enumerate}

\begin{rmrk}
  It is worth noting that the octahedron appearing in the octahedral property c) of Definition~\ref{octahedron} (which is equivalent to the braid property c') found by \cite{Fra86} in the context of homology indices of index filtrations of {\sc Morse} decompositions and exploited by him in \cite{Fra88, Fra89}) is the same as the ``octahedron diagram'' which is part of the definition of a triangulated category \cite[{IV.1}, Axiom TR4]{GM}.
\end{rmrk}

\begin{defn}[Homomorphism of graded module octahedra]
Let $(P, >)$ be a poset and let $G$ and $G'$ be graded module octahedra over $(P, >)$.
\begin{enumerate}
\item A {\em homomorphism of graded module octahedra} $\theta: G \to G'$ 
is a homomorphism of the underlying graded module triangles.
\item An {\em isomorphism of graded module octahedra} $\theta: G \to G'$ is a homomorphism of graded module octahedra where all $\theta(I)$ are isomorphisms of graded modules.
\end{enumerate}
\end{defn}

This is equivalent to the definition of homomorphisms and isomorphisms of graded module braids appearing in {\cite[Def.~1.2]{Fra88}}.

\begin{rmrk}\label{rmrk}
One notices that isomorphism of two graded module octahedra (resp.\  braids) amounts to the isomorphism of the underlying graded module triangles. Note that for each interval $I$ there is exactly one isomorphism $\theta(I):G(I)\to G'(I)$ entering in all the commutative diagrams in (\ref{Isom}) (or (\ref{isom})). This condition will be essential in Subsection \ref{obsolete}.
\end{rmrk}

\subsection{Chain complex octahedra}

\begin{defn}[Chain complex octahedron, {\cite[Def.~2.6]{Fra89}}]
Let $(P, >)$ be a poset. A {\em chain complex octahedral system} or simply a {\em chain complex octahedron} $C$ over $(P, >)$ is a chain complex triangle satisfying the extra condition:
\begin{enumerate}
\item[c)] For all $(I, J, K) \in {\mathcal I}_3(P, >)$ the following octahedron commutes and is distinguished:
 
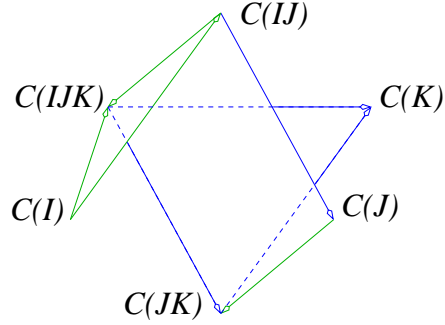
\begin{center}
      \includegraphics[scale=0.5]{octa.pstex}
    \end{center}
\end{enumerate}
\end{defn}

Since, up to quasi-isomorphism, we are talking about distinguished triangles we can complete the missing edges of the above octahedron and obtain a full octahedron as in Definition \ref{octahedron}, but now on the level of complexes (up to quasi-isomorphism). This explains again ``commutes and is distinguished''.

The two dimensional unraveling of the above octahedron is called a chain complex braid:
\[
\xymatrix@C=1.5pc@R=0.4pc{
0 \ar[dr] & & \mbox{\phantom{AA}}0\mbox{\phantom{AA}} \ar[dl] \\
& C(I) \ar[dr]^{i} \ar@/_3pc/[dd]_{i} &
& 0 \ar[dl] & \\
& &
C(IJ) \ar[dl]_{i} \ar[dr]^{p} & & 0 \ar[dl]
\\
& C(IJK) \ar[dr]^{p} \ar@/_3pc/[dd]_{p} &
&
C(J) \ar[dl]_{i} \ar[dr] & \\
& &
C(JK) \ar[dl]_{p} \ar[dr]
& & \mbox{\phantom{AA}} 0 \mbox{\phantom{AA}} \\
& C(K) \ar[dr] \ar[dl] & & 0 \\
0 & & \mbox{\phantom{AA}}0\mbox{\phantom{AA}}
}
\]

\begin{prop}[{\cite[Prop.~2.7]{Fra89}}]\label{homologychain}
Let $(P, >)$ be a poset and $C$ a chain complex octahedron (resp.\  braid) over $(P, >)$.
The homology modules $H(I)$, $I \in {\mathcal I}(P, >)$, together with
the homomorphisms $\iota(I, IJ): H(I) \to H(IJ)$ and
$\pi(IJ, J): H(IJ) \to H(J)$ induced by the chain maps
$i(I, IJ)$ and $p(IJ, J)$ form a graded module octahedron (resp.\  braid) $HC$ over $(P, >)$.
\end{prop}

A graded module octahedron (resp.\  braid) $G$ which is defined as the homology of a chain complex octahedron (resp.\  braid) is said to be {\em chain complex generated}.

\begin{prop}[{\cite[Prop.~3.4]{Fra89}}]\label{complexbraid}
Let $(P, >)$ be a poset. Further let $\{ C(p) \mid p \in P \}$ be a collection of graded modules indexed by $P$ and let $\Delta: \bigoplus_{p \in P} C(p) \to \bigoplus_{p \in P} C(p)$ be a lower triangular boundary map. Then the chain complexes $C^\Delta(I)$ defined by the $C(I)$ and $\Delta(I)$ together with the chain maps $i(I, IJ)$ and $p(IJ, J)$ form a chain complex octahedron (resp.\  braid).
\end{prop}

In other words, for a ($C$)-connection matrix for the graded module triangle $G$ to exist, $G$ must be a graded module octahedron (resp.\ braid).

\subsection{Octahedra and braids are obsolete in the definition of ($C$)-connection matrices}\label{obsolete}
When we wanted to implement the definition of connection and $C$-connection matrices in the $\Maple$ package {\tt conley}, we discovered that one can completely avoid introducing the notions of chain complex octahedra (resp.\  braids) and of graded module octahedra (resp.\  braids):\\
First notice that for $C^\Delta$ to form a chain complex braid, the matrix $\Delta$ must only satisfy the conditions of Proposition \ref{complexbraid}. By Proposition \ref{homologychain} $HC^\Delta$ is then a graded module braid. In other words, the only conditions on $\Delta$ for $HC^\Delta$ to be a graded module braid are those of Proposition \ref{complexbraid}:
\begin{enumerate}
  \item $\Delta$ is lower triangular.
  \item $\Delta$ is a boundary map.
\end{enumerate}
Furthermore Remark \ref{rmrk} says that for testing the isomorphism of $HC^\Delta$ and $G$ one only needs the isomorphism of the underlying triangles, provided that:
\begin{equation}\label{cond} \tag{$*$}
\begin{array}{l}
\mbox{For each interval $I$ there is only one isomorphism $\theta(I):HC^\Delta(I)\to G(I)$}\\
\mbox{entering in all the commutative diagrams in (\ref{Isom}) (resp.\  (\ref{isom})).}
\end{array}
\end{equation}

\subsection{{\sc Franzosa}'s existence results for $C$-connection matrices}
Now that we have collected the necessary notions we can state the main existence result for $C$-connection matrices:
\begin{axiom}[{\cite{Fra89}, Theorem 3.8\footnote{A nice coincidence.}}]
  Let $G$ be a chain complex generated graded module octahedron, $C=\{C(p) \mid  p\in P\}$ a collection of free graded modules and $\delta=\{\delta(p):C(p)\to C(p)\mid p\in P\}$ a collection of boundary maps, such that the homology of $(C(p),\delta(p))$ coincides with $G(p)$. Then there exists a $C$-connection matrix $\Delta$ for $G$ with $\Delta_{p,p}=\delta(p)$.
\end{axiom}

The converse is established by Propositions \ref{complexbraid} and \ref{homologychain}.

\begin{coro}
  If $G$ is a chain complex generated graded module octahedron and $G(p)$ is free for all $p\in P$, then there exists a connection matrix $\Delta$ for $G$.
\end{coro}
\begin{proof}
  Set $C(p)=G(p)$ with zero differential.
\end{proof}

\section{Symmetric Connection Matrices} 
\label{symm}

Let $(P,>)$ be a poset and $\Gamma$ be a subgroup of the automorphism group $\Aut(P,>)$, i.e. $q>p$ implies $\sigma q>\sigma p$ for all $\sigma\in\Gamma$. The action of $\Gamma$ on $P$ induces an action on $\mathcal{I}_n(P,>)$ for all $n\in\N$.

Let $C=\{C(p)\mid p\in P\}$ be a collection of graded modules indexed by $P$. An action of $\Gamma$ on $C$ is given by specifying for each $\sigma\in\Gamma$ and each $p\in P$ an isomorphism, which we simply denote by $\psi(\sigma):C(p)\to C(\sigma p)$, such that $\psi(\sigma\tau)=\psi(\sigma)\psi(\tau)$ for all $\sigma,\tau\in\Gamma$ and $\psi(\id_P) = \id$.

We say a linear map $\Delta:\bigoplus_{p \in P} C(p) \to \bigoplus_{p \in P} C(p)$ is $\Gamma$-symmetric, if for all $p,q\in P$ and $\sigma\in\Gamma$ the following diagram commutes:
\[
  \xymatrix@=1.5pc{
    C(q)\ar[r]^{\Delta_{q,p}}\ar[d]_{\psi(\sigma)} & C(p)\ar[d]^{\psi(\sigma)} \\
    C(\sigma q)\ar[r]^{\Delta_{\sigma q,\sigma p}} & C(\sigma p)
   }
\]

Let $C$ be a chain complex triangle. An action of $\Gamma$ on $C$ is given by specifying for each $\sigma\in\Gamma$ and each $I\in \mathcal{I}(P,>)$ an isomorphism, which we simply denote by $\psi(\sigma):C(I)\to C(\sigma I)$, such that 
\begin{enumerate}
  \item $\psi(\sigma\tau)=\psi(\sigma)\psi(\tau)$ for all $\sigma,\tau\in\Gamma$
    and $\psi(\id_P) = \id$.
  \item For all $I\in\mathcal{I}(P,>)$ and all $\sigma\in\Gamma$ the following
    diagram commutes:
    \[
      \xymatrix@=1.5pc{
        C(I)\ar[r]^{\partial(I)}\ar[d]_{\psi(\sigma)} & C(I)\ar[d]^{\psi(\sigma).} \\
        C(\sigma I)\ar[r]^{\partial(\sigma I)} & C(\sigma I)
      }
    \]
  \item For all $(I,J)\in\mathcal{I}_2(P,>)$ and all $\sigma\in\Gamma$ the following
     two diagrams commute:
     \[
       \xymatrix@=1.5pc{
         C(I)\ar[rr]^{i(I,IJ)}\ar[d]_{\psi(\sigma)} && C(IJ)\ar[d]^{\psi(\sigma),} \\
         C(\sigma I)\ar[rr]^{i(\sigma I, \sigma IJ)} && C(\sigma IJ)
       }
       \xymatrix@=1.5pc{
         C(IJ)\ar[rr]^{p(IJ,J)}\ar[d]_{\psi(\sigma)} && C(J)\ar[d]^{\psi(\sigma).} \\
         C(\sigma IJ)\ar[rr]^{p(\sigma IJ, \sigma J)} && C(\sigma J)
       }
     \]
\end{enumerate}

Let $G$ be a graded module triangle. An action of $\Gamma$ on $G$ is given by specifying for each $\sigma\in\Gamma$ and each $I\in \mathcal{I}(P,>)$ an isomorphism, which we simply denote by $\psi(\sigma):G(I)\to G(\sigma I)$, such that 
\begin{enumerate}
  \item $\psi(\sigma\tau)=\psi(\sigma)\psi(\tau)$ for all $\sigma,\tau\in\Gamma$
    and $\psi(\id_P) = \id$.
  \item For all $(I,J)\in\mathcal{I}_2(P,>)$ and all $\sigma\in\Gamma$ the following
     three diagrams commute:
     \[ \hspace{-0.8cm}
       \xymatrix@=1.5pc{
         G(I)\ar[rr]^{i(I,IJ)}\ar[d]_{\psi(\sigma)} && G(IJ)\ar[d]^{\psi(\sigma),} \\
         G(\sigma I)\ar[rr]^{i(\sigma I, \sigma IJ)} && G(\sigma IJ)
       }
       \xymatrix@=1.5pc{
         G(IJ)\ar[rr]^{p(IJ,J)}\ar[d]_{\psi(\sigma)} && G(J)\ar[d]^{\psi(\sigma),} \\
         G(\sigma IJ)\ar[rr]^{p(\sigma IJ, \sigma J)} && G(\sigma J)
       }
       \xymatrix@=1.5pc{
         G(J)\ar[rr]^{\partial(J,I)}\ar[d]_{\psi(\sigma)} && G(I)\ar[d]^{\psi(\sigma).} \\
         G(\sigma J)\ar[rr]^{\partial(\sigma J, \sigma I)} && G(\sigma I)
       }
     \]
\end{enumerate}

Note that the existence of an action of non-trivial $\Gamma\leq\Aut(P,>)$ is a non-trivial condition on $C$, resp.\  $G$.

\begin{rmrk}[Symmetric $C$-connection matrices]
  Let $\Gamma$ as above act on a collection $C=\{C(p)\mid p\in P\}$ of graded modules and let $\Delta:\bigoplus_{p \in P} C(p) \to \bigoplus_{p \in P} C(p)$ be a $\Gamma$-symmetric lower triangular boundary map. Then $\Gamma$ acts on $C^\Delta$ via $\psi(\sigma):C(I)\to C(\sigma I)$ given by the direct sum of the $\psi(\sigma):C(p)\to C(\sigma p)$, for all $p\in I$. It further acts on $HC^\Delta$ in an obvious way. Hence, for such a $\Gamma$-symmetric $\Delta$ to be a $C$-connection matrix of a graded module octahedron $G$, $\Gamma$ must act on $G$.
\end{rmrk}

More generally, if $\Gamma$ acts on a chain complex octahedron, it acts automatically on its homology graded module octahedron.

\section{$C$-Connection Matrices of a {\sc Morse} Decomposition}\label{Morse}

An inexhaustible reservoir of chain complex generated graded module octahedra are the homology graded module octahedra of a {\sc Morse} decomposition of a dynamical system. This was proved in \cite{Fra86}. See also \cite{Fra88, Fra89, Mis95, Mis} for further details, results and applications.

Let $X$ be a locally compact metric space. The object of study is
a {\em flow} $\phi: \R \times X \to X$, i.e.\ a continuous map
$\R \times X \to X$ which satisfies $\phi(0, x) = x$ and
$\phi(s, \phi(t, x)) = \phi(s+t,x)$ for all $x \in X$ and
$s, t \in \R$.

For a subset $Y\subset X$ define
\[
  \Inv(Y):=\Inv(Y,\phi):=\{x\in Y|\  \phi(\R,x)\subset Y\}\subset Y,
\]
the {\em invariant subset} of $Y$.

\subsection{Homology {\sc Conley} index}
A subset $S \subset X$ is {\em invariant} under the flow $\phi$, if $S = \Inv(S)$. It is an {\em isolated invariant set} if there exists a {\em compact} set
$Y \subset X$ (an {\em isolating neighborhood}) such that
\[
  S = \Inv(Y) \subset Y^\circ,
\]
where $Y^\circ$ denotes the interior of $Y$.

Let $M$ be an isolated invariant set \cite{Con}. A pair of compact sets $(N,L)$ with $L\subset N$ is called an {\em index pair} for $M$ if
\begin{enumerate}
  \item $\overline{N\setminus L}$ is an {\em isolating neighborhood} of $M$.
  \item $L$ is {\em positively invariant}, i.e. $\phi([0,t],x)\subset L$ for all
        $x\in L$ satisfying $\phi([0,t],x)\subset N$.
  \item $L$ is an {\em exit set} for $N$, i.e. for all $x\in N$ and all $t_1>0$ such
        that $\phi(t_1,x)\not\in N$, there exists a $t_0\in [0,t_1]$ for which
        $\phi([0,t_0],x)\subset N$ and $\phi(t_0,x)\in L$.
\end{enumerate}

Every isolated invariant set $M$ has an index pair $(N,L)$ and the {\em homology Conley index} of $M$ is defined by
\[
  CH(M):=H(N,L).
\]
$CH(M)$ is as such a graded module.

\subsection{{\sc Morse} decomposition}

For a subset $Y\subset X$ the {\em $\omega$-limit set} of $Y$ is
\[
  \omega(Y):=\bigcap_{t>0}\overline{\phi([t,\infty),Y)},
\]
while the {\em $\alpha$-limit set} of $Y$ is
\[
  \alpha(Y):=\bigcap_{t>0}\overline{\phi((-\infty,-t),Y)}.
\]

For two subsets $Y_1,Y_2\subset X$ define the {\em set of connecting orbits}
\[
  \Con(Y_1,Y_2):=\{x\in X\mid \alpha(x)\subset Y_1 \mbox{ and } \omega(x)\subset Y_2\}.
\]

Let $S$ be an isolated invariant set and $(P,>)$ be a poset. A finite collection
\[ 
  \mathcal{M}(S)=\{M(p)\mid p\in P\}
\]
of disjoint isolated invariant subsets of $S$ is called a {\em {\sc Morse} decomposition} if there exists a strict partial order $>$ on $P$, such that for every $x\in S\setminus \bigcup_{p\in P}M(p)$ there exists $p,q\in P$, such that $q>p$ and $x\in\Con(M(q),M(p))$.

The sets $M(p)$ are called {\em {\sc Morse} sets}. A partial order on $P$ satisfying this property is said to be {\em admissible}.

There is a partial order $>_\phi$ induced by the flow, {\em generated} by the relations $q>_\phi p$ whenever $\Con(M(q),M(p))\neq\emptyset$. This so called flow-induced order is a subset of every admissible order, and in this sense minimal. Normally this order is not known and one is content with a coarser order. If, for example, an (energy) function $E$ is known with $E(x) > E(\phi(t,x))$ for all $t>0$, then defining the partial order $>_E$ by
\[
  q>_E p, \mbox{ iff $E(y) > E(x)$ for all $y\in M(q)$ and $x\in M(p)$},
\] 
yields an admissible order.

For an interval $I$ define the {\sc Morse} set
\[
  M(I):=\bigcup_{p\in I}M(p)\cup\bigcup_{p,q\in I}\Con(M(q),M(p)).
\]
$M(I)$ is again an isolated invariant set. If $(I,J)\in\mathcal{I}_2(P,>)$, then $(M(I),M(J))$ is an attractor-repeller pair in $M(IJ)$.

\begin{rmrk}\label{indexfiltration}
In dynamical system theory {\sc Franzosa} discovered in \cite{Fra86} a natural chain complex octahedron, the so called chain complex octahedron $C^\NN$ defined by a so called {\em index filtration} $\NN$ of a {\sc Morse} decomposition. Its homology $HC^\NN$ is by Proposition~\ref{homologychain} a (chain complex generated) graded module octahedron:

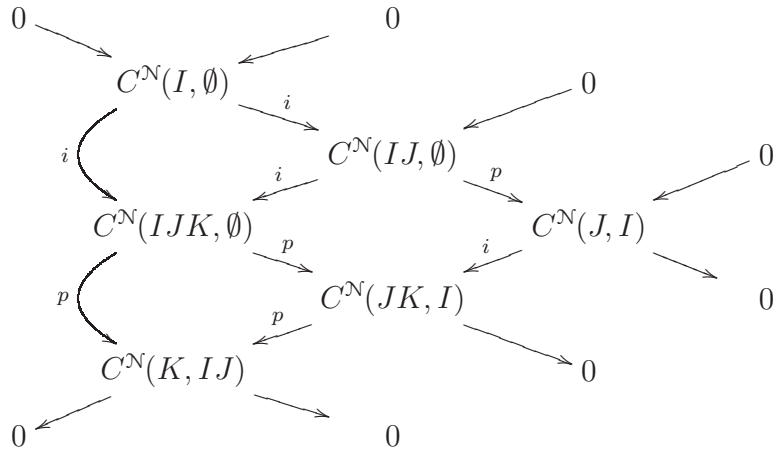
\[
\xymatrix@C=1.5pc@R=0.6pc{
0 \ar[dr] & & \mbox{\phantom{AA}}0\mbox{\phantom{AA}} \ar[dl] \\
& C^\NN(I,\emptyset) \ar[dr]^{i} \ar@/_3pc/[dd]_{i} &
& 0 \ar[dl] & \\
& &
C^\NN(IJ,\emptyset) \ar[dl]_{i} \ar[dr]^{p} & & 0 \ar[dl]
\\
& C^\NN(IJK,\emptyset) \ar[dr]^{p} \ar@/_3pc/[dd]_{p} &
&
C^\NN(J,I) \ar[dl]_{i} \ar[dr] & \\
& &
C^\NN(JK,I) \ar[dl]_{p} \ar[dr]
& & \mbox{\phantom{AA}} 0 \mbox{\phantom{AA}} \\
& C^\NN(K,IJ) \ar[dr] \ar[dl] & & 0 \\
0 & & \mbox{\phantom{AA}}0\mbox{\phantom{AA}}
}
\]
where $C^\NN(A,B):=C(\NN(QAB),\NN(QB))$ for $A,B$ one of the above pairs of intervals and $Q$ is a certain interval defined by $(I,J,K)\in\mathcal{I}_3(P,>)$. For the reasons mentioned in Remark \ref{rmrk2} below we won't go into the details of this construction. The only thing one has to know is that $C(\NN(QAB),\NN(QB))$ is the relative singular (or simplicial, ...) chain complex\footnote{In talking about the chain complex braid we use the fact that the quasi-isomorphism type of $C^\NN(A,B)$ neither depends on $B$ nor on $Q$.} of the pair of compact topological spaces $(\NN(QAB),\NN(QB))$, which is an index pair for the {\sc Morse} set $M(A)$ (independent of $B$ and $Q$!), i.e. its relative homology being the homological {\sc Conley} index $CH(M(A))$ of $M(A)$.
\end{rmrk}

\begin{defn}[($C$)-Connection matrix of a {\sc Morse} decomposition]
  Given a {\sc Morse} decomposition with an index filtration $\NN$. A {\em ($C$)-connection matrix} of the {\sc Morse} decomposition is a ($C$)-connection matrix of $HC^\NN$. One can show that the definition is independent of the index filtration.
\end{defn}

At this point we want to emphasize that in order to comply with Definition \ref{conn}, one needs to correctly identify the $CH(M(A))$ arising as the homology of the different, but quasi-isomorphic complexes coming from the topological setup.

Remark \ref{indexfiltration} states that $G=HC^\NN$ is automatically a chain complex generated graded module octahedron. In other words, in the context of ($C$)-connection matrices of {\sc Morse} decompositions, one never needs to check the octahedron (resp.\  braid) condition. Cf. \cite[p.~574, below Theorem 3.8]{Fra86}.

It remains to mention why the theory of ($C$)-connection matrices is relevant in the context of {\sc Morse} decompositions. There is a link between existence of connecting orbits and non-triviality of entries of connection matrices. For details see for example \cite[4.1.1]{Mis}.

\section{The $\Maple$ Package $\conley$}
\label{package}

Although the theory of $C$-connection matrices is a purely algebraic theory, we've chosen $\conley$ as the name of our $\Maple$ package in honor of {\sc Charles Conley}, since his index theory \cite{Con} is a natural field for applying these ideas. As mentioned above, the case $G=HC^\NN$ is of special interest to applications.

In {\sc Conley} index theory one uses as much topological data as possible given by (the index filtration of) the {\sc Morse} decomposition to restrict the possible connection matrices $\Delta$ for $HC^\NN$.

\begin{rmrk}[Even less is used in practice]\label{rmrk2}
Typically, the only data that are given are the isomorphism types of the homology {\sc Conley} indices $CH(M(I))$ for {\em some} $I\in \mathcal{I}(P,>)$. Nevertheless, the minimal data required remains the isomorphism type of the homology {\sc Conley} indices for {\em all} one element intervals.

As mentioned at the end of Remark \ref{indexfiltration}, one needs, in order to apply the definition of a connection matrix, to identify for each interval $I$ the homologies of different index pairs for $M(I)$ inside the index filtration, which are all isomorphic to $CH(M(I))$. Since this is an extremely demanding task, and up to our knowledge not addressed yet, one cannot make use of condition (\ref{cond}).

To be able to make use of the commutativity of the squares in (\ref{Isom}) (resp.\  (\ref{isom})) one needs the induced maps coming from the index triples for each pair of adjacent intervals, a data which is rarely provided. 

Summing up, in the majority of non-trivial examples, the data of the index filtration for the {\sc Morse} decomposition are up to our knowledge not known completely. In particular one can neither make use of condition (\ref{cond}) nor even of the commutativity of the squares in (\ref{Isom}) (resp.\  (\ref{isom})) to impose further restrictions on $\Delta$. The only restriction left on $\Delta$ is the abstract isomorphism between $HC^\Delta(I)$ and $CH(M(I))$, in case the isomorphism type of the latter is known. Due to this lack of topological data the package {\tt conley} only imposes the last mentioned restriction on $\Delta$. If at some point in the future such data are provided, {\tt conley} can easily be extended to make use of them.
\end{rmrk}

\section{Examples}
\label{beispiele}

In each of the four examples in this section the set of all connection matrices which are compatible with the given data of the flow under consideration are computed using the $\Maple$ package $\conley$. As $\conley$ is based on $\homalg$, the computations on the level of homological algebra are carried out by $\homalg$. Since in the following examples we compute with modules over a principal ideal ring, the so-called ring package $\mathtt{PIR}$ is used to do the ring arithmetics for $\homalg$ \cite{BR}.

For the $\conley$ package the following conventions to provide graded modules, e.g.\ homology modules, are used:
\begin{itemize}
\item The most elaborate notation is to provide the graded components as list of presentations in the $\homalg$ format.
\item If all graded components are free modules then a list of ranks suffices.
\item If the graded module is of rank $1$ and concentrated in one degree then this degree is sufficient as input. In this case, we call this degree the
{\em index} of the graded module.
\end{itemize}

\subsection{A flow with $V_4\times C_2\cong C_2^3$-symmetry}

\PUSH{V4xC2.tex}%
\DefineParaStyle{Left Justified Maple Output}
\DefineParaStyle{Maple Output}
\DefineParaStyle{_pstyle1}
\DefineParaStyle{_pstyle8}
\DefineCharStyle{2D Output}
\DefineCharStyle{LaTeX}
\begin{maplegroup}
We consider a $C_2^3$-symmetric flow on the $2$-sphere with six
hyperbolic equilibria. They are enumerated as shown in the following
figure: 
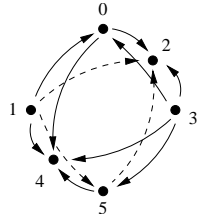
\begin{center}  \includegraphics[scale=0.6]{V4xC2.pstex}
\end{center}

\end{maplegroup}
\begin{maplegroup}
\begin{_pstyle1}
\begin{mapleinput}
\mapleinline{active}{1d}{restart;}{%
}
\end{mapleinput}
\end{_pstyle1}

\end{maplegroup}
\begin{maplegroup}
\begin{_pstyle1}
\begin{mapleinput}
\mapleinline{active}{1d}{with(conley): with(PIR): with(homalg):}{%
}
\end{mapleinput}
\end{_pstyle1}

\end{maplegroup}
\begin{maplegroup}
\begin{_pstyle1}
\begin{mapleinput}
\mapleinline{active}{1d}{`homalg/default`:=`PIR/homalg`:}{%
}
\end{mapleinput}
\end{_pstyle1}

\end{maplegroup}
\begin{maplegroup}
We want to compute over $\Z/2 \Z$:

\end{maplegroup}
\begin{maplegroup}
\begin{_pstyle1}
\begin{mapleinput}
\mapleinline{active}{1d}{var:=[2];}{%
}
\end{mapleinput}
\end{_pstyle1}

\mapleresult
\begin{maplelatex}
\mapleinline{inert}{2d}{var := [2];}{%
\[
\mathit{var} := [2]
\]
}
\end{maplelatex}

\end{maplegroup}
\begin{maplegroup}
\begin{_pstyle1}
\begin{mapleinput}
\mapleinline{active}{1d}{Pvar(var);}{%
}
\end{mapleinput}
\end{_pstyle1}

\mapleresult
\begin{maplelatex}
\mapleinline{inert}{2d}{["Z/pZ", 2];}{%
\[
[\mbox{``Z/pZ''}, \,2]
\]
}
\end{maplelatex}

\end{maplegroup}
\begin{maplegroup}
The set $P$ of six equilibria of the flow is defined:

\end{maplegroup}
\begin{maplegroup}
\begin{mapleinput}
\mapleinline{active}{1d}{P:=[0,1,2,3,4,5];}{%
}
\end{mapleinput}

\mapleresult
\begin{maplelatex}
\mapleinline{inert}{2d}{P := [0, 1, 2, 3, 4, 5];}{%
\[
P := [0, \,1, \,2, \,3, \,4, \,5]
\]
}
\end{maplelatex}

\end{maplegroup}
\begin{maplegroup}
The symmetry group $\Gamma := V_4 \times C_2 \cong C_2^3$ is generated
by the following permutations:

\end{maplegroup}
\begin{maplegroup}
\begin{mapleinput}
\mapleinline{active}{1d}{V4C2:=[ [[1,3]], [[2,4]], [[0,5]] ];}{%
}
\end{mapleinput}

\mapleresult
\begin{maplelatex}
\mapleinline{inert}{2d}{V4C2 := [[[1, 3]], [[2, 4]], [[0, 5]]];}{%
\[
\mathit{V4C2} := [[[1, \,3]], \,[[2, \,4]], \,[[0, \,5]]]
\]
}
\end{maplelatex}

\end{maplegroup}
\begin{maplegroup}
The generating set of relations for $>$ is given by:

\end{maplegroup}
\begin{maplegroup}
\begin{_pstyle1}
\begin{mapleinput}
\mapleinline{active}{1d}{rel_:=[[1,0],[0,2]];}{%
}
\end{mapleinput}
\end{_pstyle1}

\mapleresult
\begin{maplelatex}
\mapleinline{inert}{2d}{rel_ := [[1, 0], [0, 2]];}{%
\[
\mathit{rel\_} := [[1, \,0], \,[0, \,2]]
\]
}
\end{maplelatex}

\end{maplegroup}
\begin{maplegroup}
The $\Gamma$-symmetry gives more relations:

\end{maplegroup}
\begin{maplegroup}
\begin{mapleinput}
\mapleinline{active}{1d}{Orbits(V4C2,rel_):
rel:=map(op,\%);}{%
}
\end{mapleinput}

\mapleresult
\begin{maplelatex}
\mapleinline{inert}{2d}{rel := [[1, 0], [3, 0], [1, 5], [3, 5], [0, 2], [0, 4], [5, 2], [5,
4]];}{%
\[
\mathit{rel} := [[1, \,0], \,[3, \,0], \,[1, \,5], \,[3, \,5], \,
[0, \,2], \,[0, \,4], \,[5, \,2], \,[5, \,4]]
\]
}
\end{maplelatex}

\end{maplegroup}
\begin{maplegroup}
\begin{_pstyle1}
\begin{mapleinput}
\mapleinline{active}{1d}{ord:=GeneratePartialOrder(P,rel);}{%
}
\end{mapleinput}
\end{_pstyle1}

\mapleresult
\begin{maplelatex}
\mapleinline{inert}{2d}{ord := [[0, 2], [0, 4], [1, 0], [1, 2], [1, 4], [1, 5], [3, 0], [3,
2], [3, 4], [3, 5], [5, 2], [5, 4]];}{%
\[
\mathit{ord} := [[0, \,2], \,[0, \,4], \,[1, \,0], \,[1, \,2], \,
[1, \,4], \,[1, \,5], \,[3, \,0], \,[3, \,2], \,[3, \,4], \,[3, 
\,5], \,[5, \,2], \,[5, \,4]]
\]
}
\end{maplelatex}

\end{maplegroup}
\begin{maplegroup}
\begin{_pstyle1}
\begin{mapleinput}
\mapleinline{active}{1d}{I1:=GenerateIntervals(P,rel);}{%
}
\end{mapleinput}
\end{_pstyle1}

\mapleresult
\begin{maplelatex}
\mapleinline{inert}{2d}{I1 := [[], [2], [1], [5], [0], [3], [4], [1, 3], [2, 4], [0, 5], [0,
2], [1, 5], [3, 5], [0, 4], [0, 1], [0, 3], [2, 5], [4, 5], [0, 4, 5],
[2, 4, 5], [0, 1, 5], [0, 2, 5], [0, 3, 5], [1, 3, 5], [0, 1, 3], [0,
2, 4], [0, 3, 4, 5], [0, 1, 4, 5], [0, 2, 4, 5], [0, 1, 2, 5], [0, 1,
3, 5], [0, 2, 3, 5], [0, 1, 3, 4, 5], [0, 2, 3, 4, 5], [0, 1, 2, 4,
5], [0, 1, 2, 3, 5], [0, 1, 2, 3, 4, 5]];}{%
\maplemultiline{
\mathit{I1} := [[], \,[2], \,[1], \,[5], \,[0], \,[3], \,[4], \,[
1, \,3], \,[2, \,4], \,[0, \,5], \,[0, \,2], \,[1, \,5], \,[3, \,
5], \,[0, \,4], \,[0, \,1],  \\
[0, \,3], \,[2, \,5], \,[4, \,5], \,[0, \,4, \,5], \,[2, \,4, \,5
], \,[0, \,1, \,5], \,[0, \,2, \,5], \,[0, \,3, \,5], \,[1, \,3, 
\,5],  \\
[0, \,1, \,3], \,[0, \,2, \,4], \,[0, \,3, \,4, \,5], \,[0, \,1, 
\,4, \,5], \,[0, \,2, \,4, \,5], \,[0, \,1, \,2, \,5], \,[0, \,1
, \,3, \,5],  \\
[0, \,2, \,3, \,5], \,[0, \,1, \,3, \,4, \,5], \,[0, \,2, \,3, \,
4, \,5], \,[0, \,1, \,2, \,4, \,5], \,[0, \,1, \,2, \,3, \,5], 
 \\
[0, \,1, \,2, \,3, \,4, \,5]] }
}
\end{maplelatex}

\end{maplegroup}
\begin{maplegroup}
\begin{_pstyle1}
\begin{mapleinput}
\mapleinline{active}{1d}{#I123:=GenerateAdjacentIntervals(P,rel);}{%
}
\end{mapleinput}
\end{_pstyle1}

\end{maplegroup}
\begin{maplegroup}
The $\Gamma$-orbits on $> \, \subset P \times P$:

\end{maplegroup}
\begin{maplegroup}
\begin{_pstyle1}
\begin{mapleinput}
\mapleinline{active}{1d}{orbits:=Orbits(V4C2,ord);}{%
}
\end{mapleinput}
\end{_pstyle1}

\mapleresult
\begin{maplelatex}
\mapleinline{inert}{2d}{orbits := [[[0, 2], [0, 4], [5, 2], [5, 4]], [[1, 0], [3, 0], [1, 5],
[3, 5]], [[1, 2], [3, 2], [1, 4], [3, 4]]];}{%
\maplemultiline{
\mathit{orbits} := [[[0, \,2], \,[0, \,4], \,[5, \,2], \,[5, \,4]
], \,[[1, \,0], \,[3, \,0], \,[1, \,5], \,[3, \,5]],  \\
[[1, \,2], \,[3, \,2], \,[1, \,4], \,[3, \,4]]] }
}
\end{maplelatex}

\end{maplegroup}
\begin{maplegroup}
Now we have to provide the homology {\sc Conley} indices for the six
equilibria:

\end{maplegroup}
\begin{maplegroup}
\begin{_pstyle1}
\begin{mapleinput}
\mapleinline{active}{1d}{CHp:=[1,2,0,2,0,1];}{%
}
\end{mapleinput}
\end{_pstyle1}

\mapleresult
\begin{maplelatex}
\mapleinline{inert}{2d}{CHp := [1, 2, 0, 2, 0, 1];}{%
\[
\mathit{CHp} := [1, \,2, \,0, \,2, \,0, \,1]
\]
}
\end{maplelatex}

\end{maplegroup}
\begin{maplegroup}
The homology of the sphere is added as homology {\sc Conley} index of
$M(P)$:

\end{maplegroup}
\begin{maplegroup}
\begin{_pstyle1}
\begin{mapleinput}
\mapleinline{active}{1d}{CH:=[op(zip((x,y)->x=y,P,CHp)),P=[1,0,1]];}{%
}
\end{mapleinput}
\end{_pstyle1}

\mapleresult
\begin{maplelatex}
\mapleinline{inert}{2d}{CH := [0 = 1, 1 = 2, 2 = 0, 3 = 2, 4 = 0, 5 = 1, [0, 1, 2, 3, 4, 5] =
[1, 0, 1]];}{%
\[
\mathit{CH} := [0=1, \,1=2, \,2=0, \,3=2, \,4=0, \,5=1, \,[0, \,1
, \,2, \,3, \,4, \,5]=[1, \,0, \,1]]
\]
}
\end{maplelatex}

\end{maplegroup}
\begin{maplegroup}
The following computation of all connection matrices takes into
account the $\Gamma$-symmetry. The optional string
$\mathtt{"full\_notation"}$ means that for the matrix
$\Delta_i:\bigoplus_{p\in P} CH_i(p)\to \bigoplus_{p\in P}
CH_{i-1}(p)$ a zero row is reserved for every trivial $CH_i(p)$ and a
zero column is reserved for every trivial $CH_{i-1}(p)$. In the
hyperbolic case this leads to $|P|\times|P|$-square matrices:

\end{maplegroup}
\begin{maplegroup}
\begin{_pstyle1}
\begin{mapleinput}
\mapleinline{active}{1d}{Con:=ConnectionMatrices(P,rel,CH,var,"Symmetry"=V4C2,"full_notation")
;}{%
}
\end{mapleinput}
\end{_pstyle1}

\mapleresult
\begin{maplelatex}
\mapleinline{inert}{2d}{Con := [[matrix([[0, 0, 1, 0, 1, 0], [0, 0, 0, 0, 0, 0], [0, 0, 0, 0,
0, 0], [0, 0, 0, 0, 0, 0], [0, 0, 0, 0, 0, 0], [0, 0, 1, 0, 1, 0]]),
matrix([[0, 0, 0, 0, 0, 0], [1, 0, 0, 0, 0, 1], [0, 0, 0, 0, 0, 0],
[1, 0, 0, 0, 0, 1], [0, 0, 0, 0, 0, 0], [0, 0, 0, 0, 0, 0]])]];}{%
\[
\mathit{Con} :=  \left[  \!  \left[  \!  \left[ 
{\begin{array}{rrrrrr}
0 & 0 & 1 & 0 & 1 & 0 \\
0 & 0 & 0 & 0 & 0 & 0 \\
0 & 0 & 0 & 0 & 0 & 0 \\
0 & 0 & 0 & 0 & 0 & 0 \\
0 & 0 & 0 & 0 & 0 & 0 \\
0 & 0 & 1 & 0 & 1 & 0
\end{array}}
 \right] , \, \left[ 
{\begin{array}{rrrrrr}
0 & 0 & 0 & 0 & 0 & 0 \\
1 & 0 & 0 & 0 & 0 & 1 \\
0 & 0 & 0 & 0 & 0 & 0 \\
1 & 0 & 0 & 0 & 0 & 1 \\
0 & 0 & 0 & 0 & 0 & 0 \\
0 & 0 & 0 & 0 & 0 & 0
\end{array}}
 \right]  \!  \right]  \!  \right] 
\]
}
\end{maplelatex}

\end{maplegroup}
\begin{maplegroup}
\begin{mapleinput}
\mapleinline{active}{1d}{nops(Con);}{%
}
\end{mapleinput}

\mapleresult
\begin{maplelatex}
\mapleinline{inert}{2d}{1;}{%
\[
1
\]
}
\end{maplelatex}

\end{maplegroup}
\begin{maplegroup}
We conclude that there is exactly one connection matrix which is
compatible with the given data. It is given by the list containing the
matrices which represent the boundary maps from $\bigoplus_{p \in P}
CH_1(p)$ to $\bigoplus_{p \in P} CH_0(p)$ resp.\  from $\bigoplus_{p
\in P} CH_2(p)$ to $\bigoplus_{p \in P} CH_1(p)$. Finally, we check
that connecting orbits from $1$ to $0$ resp.\ from $0$ to $2$ exist:

\end{maplegroup}
\begin{maplegroup}
\begin{mapleinput}
\mapleinline{active}{1d}{map(ConnectionSubmatrix,Con,P,[1,0],CHp,var);}{%
}
\end{mapleinput}

\mapleresult
\begin{maplelatex}
\mapleinline{inert}{2d}{[[matrix([[0]]), matrix([[1]])]];}{%
\[
[[ \left[ 
{\begin{array}{r}
0
\end{array}}
 \right] , \, \left[ 
{\begin{array}{r}
1
\end{array}}
 \right] ]]
\]
}
\end{maplelatex}

\end{maplegroup}
\begin{maplegroup}
\begin{mapleinput}
\mapleinline{active}{1d}{map(ConnectionSubmatrix,Con,P,[0,2],CHp,var);}{%
}
\end{mapleinput}

\mapleresult
\begin{maplelatex}
\mapleinline{inert}{2d}{[[matrix([[1]]), matrix([[0]])]];}{%
\[
[[ \left[ 
{\begin{array}{r}
1
\end{array}}
 \right] , \, \left[ 
{\begin{array}{r}
0
\end{array}}
 \right] ]]
\]
}
\end{maplelatex}

\end{maplegroup}
\begin{maplegroup}
Since the format of the connection matrices may change in the future,
the user is recommended to use the above way of accessing the entries
of the connection matrices.

\end{maplegroup}
%
\POP

\subsection{A flow with $D_6\times C_2\cong D_{12}$-symmetry}\label{D_12}

\PUSH{D6xC2.tex}%
\DefineParaStyle{Left Justified Maple Output}
\DefineParaStyle{Maple Output}
\DefineParaStyle{_pstyle1}
\DefineParaStyle{_pstyle8}
\DefineCharStyle{2D Output}
\DefineCharStyle{LaTeX}
\begin{maplegroup}
Here we use $\conley$ to compute the connection matrices of a
$D_{12}$-symmetric flow on the $2$-sphere with six hyperbolic and two
non-hyperbolic equilibria. They are enumerated as shown in the
following figure: 
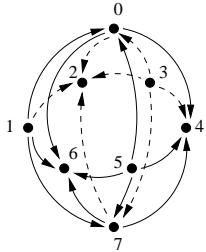
\begin{center}  
\includegraphics[scale=0.6]{D6xC2.pstex} \end{center} The following
two figures give a simple view of the upper resp.\  lower hemisphere:
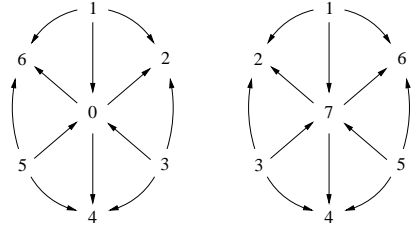
\begin{center}  \includegraphics[scale=0.6]{D6xC2_0.pstex} \qquad
\includegraphics[scale=0.6]{D6xC2_7.pstex} \end{center} Around $0$
resp.\  $7$ an isolating neighborhood looks like: 
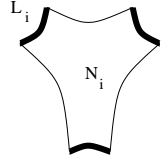
\begin{center}  
\includegraphics[scale=0.6]{LN.pstex} \end{center} $(N_i,L_i)$ is an
index pair for $i=0,7$, where $L_i$ is the disjoint union of the three
bold border regions of the isolating neighborhood $N_i$ of $i$. The
index pair $(N_i,L_i)$ is homotopy equivalent to the simplicial
complex: 
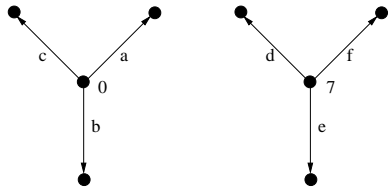
\begin{center} 
\includegraphics[scale=0.6]{D6xC2_triang0.pstex} \qquad
\includegraphics[scale=0.6]{D6xC2_triang7.pstex} \end{center}

\end{maplegroup}
\begin{maplegroup}
\begin{_pstyle1}
\begin{mapleinput}
\mapleinline{active}{1d}{restart;}{%
}
\end{mapleinput}
\end{_pstyle1}

\end{maplegroup}
\begin{maplegroup}
\begin{_pstyle1}
\begin{mapleinput}
\mapleinline{active}{1d}{with(conley): with(PIR): with(homalg):}{%
}
\end{mapleinput}
\end{_pstyle1}

\end{maplegroup}
\begin{maplegroup}
\begin{_pstyle1}
\begin{mapleinput}
\mapleinline{active}{1d}{`homalg/default`:=`PIR/homalg`:}{%
}
\end{mapleinput}
\end{_pstyle1}

\end{maplegroup}
\begin{maplegroup}
We set up the field of coefficients to be $\F_2 := \Z / 2 \Z$:

\end{maplegroup}
\begin{maplegroup}
\begin{_pstyle1}
\begin{mapleinput}
\mapleinline{active}{1d}{var:=[2];}{%
}
\end{mapleinput}
\end{_pstyle1}

\mapleresult
\begin{maplelatex}
\mapleinline{inert}{2d}{var := [2];}{%
\[
\mathit{var} := [2]
\]
}
\end{maplelatex}

\end{maplegroup}
\begin{maplegroup}
\begin{_pstyle1}
\begin{mapleinput}
\mapleinline{active}{1d}{Pvar(var);}{%
}
\end{mapleinput}
\end{_pstyle1}

\mapleresult
\begin{maplelatex}
\mapleinline{inert}{2d}{["Z/pZ", 2];}{%
\[
[\mbox{``Z/pZ''}, \,2]
\]
}
\end{maplelatex}

\end{maplegroup}
\begin{maplegroup}
We define the set $P$ of equilibria and the generators for $\Gamma :=
D_{12}$:

\end{maplegroup}
\begin{maplegroup}
\begin{mapleinput}
\mapleinline{active}{1d}{P:=[0,1,2,3,4,5,6,7];}{%
}
\end{mapleinput}

\mapleresult
\begin{maplelatex}
\mapleinline{inert}{2d}{P := [0, 1, 2, 3, 4, 5, 6, 7];}{%
\[
P := [0, \,1, \,2, \,3, \,4, \,5, \,6, \,7]
\]
}
\end{maplelatex}

\end{maplegroup}
\begin{maplegroup}
\begin{mapleinput}
\mapleinline{active}{1d}{D6C2:=[ [[1,3,5],[2,4,6]], [[3,5],[2,6]] ,[[0,7]] ];}{%
}
\end{mapleinput}

\mapleresult
\begin{maplelatex}
\mapleinline{inert}{2d}{D6C2 := [[[1, 3, 5], [2, 4, 6]], [[3, 5], [2, 6]], [[0, 7]]];}{%
\[
\mathit{D6C2} := [[[1, \,3, \,5], \,[2, \,4, \,6]], \,[[3, \,5], 
\,[2, \,6]], \,[[0, \,7]]]
\]
}
\end{maplelatex}

\end{maplegroup}
\begin{maplegroup}
The action of the generating permutations $(1, 3, 5)(2, 4, 6)$, $(3,
5)(2, 6)$ and $(0, 7)$ on $\bigoplus_{p \in P} CH_i(p)$ for $i = 0, 1,
2$ is represented by the following matrices. More precisely, for each
generator we give an equation whose left hand side is the permutation
$\pi$ of the set $P$ and whose right hand side is the corresponding
list $[M_0(\pi), M_1(\pi), M_2(\pi)]$ of matrices representing the
action of $\pi$ on the respective $\bigoplus_{p \in P} CH_i(p)$, $i =
0, 1, 2$. The dimensions of the matrices are coherent with the so
called
standard notation ($\mathtt{"std\_notation"}$) explained below:

\end{maplegroup}
\begin{maplegroup}
\begin{mapleinput}
\mapleinline{active}{1d}{id:=linalg[diag](1$3);
rho:=matrix([[0,1,0],[0,0,1],[1,0,0]]);
rep:=[ [rho,matrix([[0,1,0,0],[1,1,0,0],[0,0,0,1],[0,0,1,1]]),rho],
[matrix([[0,0,1],[0,1,0],[1,0,0]]),matrix([[1,1,0,0],[0,1,0,0],[0,0,1,
1],[0,0,0,1]]),matrix([[1,0,0],[0,0,1],[0,1,0]])],
[id,matrix([[0,0,1,0],[0,0,0,1],[1,0,0,0],[0,1,0,0]]),id]]:}{%
}
\end{mapleinput}

\mapleresult
\begin{maplelatex}
\mapleinline{inert}{2d}{id := matrix([[1, 0, 0], [0, 1, 0], [0, 0, 1]]);}{%
\[
\mathit{id} :=  \left[ 
{\begin{array}{rrr}
1 & 0 & 0 \\
0 & 1 & 0 \\
0 & 0 & 1
\end{array}}
 \right] 
\]
}
\end{maplelatex}

\begin{maplelatex}
\mapleinline{inert}{2d}{rho := matrix([[0, 1, 0], [0, 0, 1], [1, 0, 0]]);}{%
\[
\rho  :=  \left[ 
{\begin{array}{rrr}
0 & 1 & 0 \\
0 & 0 & 1 \\
1 & 0 & 0
\end{array}}
 \right] 
\]
}
\end{maplelatex}

\end{maplegroup}
\begin{maplegroup}
\begin{mapleinput}
\mapleinline{active}{1d}{D6C2_rep:=zip((a,b)->a=b,D6C2,rep);}{%
}
\end{mapleinput}

\mapleresult
\begin{maplelatex}
\mapleinline{inert}{2d}{D6C2_rep := [[[1, 3, 5], [2, 4, 6]] = [rho, matrix([[0, 1, 0, 0], [1,
1, 0, 0], [0, 0, 0, 1], [0, 0, 1, 1]]), rho], [[3, 5], [2, 6]] =
[matrix([[0, 0, 1], [0, 1, 0], [1, 0, 0]]), matrix([[1, 1, 0, 0], [0,
1, 0, 0], [0, 0, 1, 1], [0, 0, 0, 1]]), matrix([[1, 0, 0], [0, 0, 1],
[0, 1, 0]])], [[0, 7]] = [id, matrix([[0, 0, 1, 0], [0, 0, 0, 1], [1,
0, 0, 0], [0, 1, 0, 0]]), id]];}{%
\maplemultiline{
\mathit{D6C2\_rep} :=  \left[ {\vrule 
height1.47em width0em depth1.47em} \right. \!  \! [[1, \,3, \,5]
, \,[2, \,4, \,6]]= \left[  \! \rho , \, \left[ 
{\begin{array}{rrrr}
0 & 1 & 0 & 0 \\
1 & 1 & 0 & 0 \\
0 & 0 & 0 & 1 \\
0 & 0 & 1 & 1
\end{array}}
 \right] , \,\rho  \!  \right] ,  \\
[[3, \,5], \,[2, \,6]]= \left[  \!  \left[ 
{\begin{array}{rrr}
0 & 0 & 1 \\
0 & 1 & 0 \\
1 & 0 & 0
\end{array}}
 \right] , \, \left[ 
{\begin{array}{rrrr}
1 & 1 & 0 & 0 \\
0 & 1 & 0 & 0 \\
0 & 0 & 1 & 1 \\
0 & 0 & 0 & 1
\end{array}}
 \right] , \, \left[ 
{\begin{array}{rrr}
1 & 0 & 0 \\
0 & 0 & 1 \\
0 & 1 & 0
\end{array}}
 \right]  \!  \right] , \\
\,[[0, \,7]]= \left[  \! \mathit{id}, \,
 \left[ 
{\begin{array}{rrrr}
0 & 0 & 1 & 0 \\
0 & 0 & 0 & 1 \\
1 & 0 & 0 & 0 \\
0 & 1 & 0 & 0
\end{array}}
 \right] , \,\mathit{id} \!  \right]  \! \! \left. {\vrule 
height1.47em width0em depth1.47em} \right]  }
}
\end{maplelatex}

\end{maplegroup}
\begin{maplegroup}
A simple calculation shows that

\[ CH_1(0):=H_1(N_0,L_0)\cong \F_2 a\oplus \F_2 b \oplus \F_2 c/ \F_2
(a+b+c) = \langle\overline{a},\overline{b}\rangle_{\F_2}\cong
\F_2^{1\times 2} \] resp.\  
\[ CH_1(7):=H_1(N_7,L_7)\cong \F_2 d\oplus \F_2 e \oplus \F_2 f/ \F_2
(d+e+f) = \langle\overline{d},\overline{e}\rangle_{\F_2}\cong
\F_2^{1\times 2}. \] The $\F_2$-basis of $\bigoplus_{p \in P}
CH_1(p)=CH_1(0)\oplus CH_1(7)\cong \F_2^{1\times 2}\oplus
\F_2^{1\times 2}\cong \F_2^{1\times 4}$ which was chosen to write down
the above matrices is
$(\overline{a},\overline{b},\overline{d},\overline{e})$.

\end{maplegroup}
\begin{maplegroup}
Next we define the generating set of relations for $>$:

\end{maplegroup}
\begin{maplegroup}
\begin{_pstyle1}
\begin{mapleinput}
\mapleinline{active}{1d}{rel_:=[[1,0],[0,2]];}{%
}
\end{mapleinput}
\end{_pstyle1}

\mapleresult
\begin{maplelatex}
\mapleinline{inert}{2d}{rel_ := [[1, 0], [0, 2]];}{%
\[
\mathit{rel\_} := [[1, \,0], \,[0, \,2]]
\]
}
\end{maplelatex}

\end{maplegroup}
\begin{maplegroup}
The orbits of these relations under the given symmetry group are
obtained as follows:

\end{maplegroup}
\begin{maplegroup}
\begin{mapleinput}
\mapleinline{active}{1d}{Orbits(D6C2,rel_):
rel:=map(op,\%);}{%
}
\end{mapleinput}

\mapleresult
\begin{maplelatex}
\mapleinline{inert}{2d}{rel := [[1, 0], [3, 0], [1, 7], [5, 0], [3, 7], [5, 7], [0, 2], [0,
4], [0, 6], [7, 2], [7, 4], [7, 6]];}{%
\[
\mathit{rel} := [[1, \,0], \,[3, \,0], \,[1, \,7], \,[5, \,0], \,
[3, \,7], \,[5, \,7], \,[0, \,2], \,[0, \,4], \,[0, \,6], \,[7, 
\,2], \,[7, \,4], \,[7, \,6]]
\]
}
\end{maplelatex}

\end{maplegroup}
\begin{maplegroup}
\begin{_pstyle1}
\begin{mapleinput}
\mapleinline{active}{1d}{ord:=GeneratePartialOrder(P,rel);}{%
}
\end{mapleinput}
\end{_pstyle1}

\mapleresult
\begin{maplelatex}
\mapleinline{inert}{2d}{ord := [[0, 2], [0, 4], [0, 6], [1, 0], [1, 2], [1, 4], [1, 6], [1,
7], [3, 0], [3, 2], [3, 4], [3, 6], [3, 7], [5, 0], [5, 2], [5, 4],
[5, 6], [5, 7], [7, 2], [7, 4], [7, 6]];}{%
\maplemultiline{
\mathit{ord} := [[0, \,2], \,[0, \,4], \,[0, \,6], \,[1, \,0], \,
[1, \,2], \,[1, \,4], \,[1, \,6], \,[1, \,7], \,[3, \,0], \,[3, 
\,2], \,[3, \,4], \,[3, \,6],  \\
[3, \,7], \,[5, \,0], \,[5, \,2], \,[5, \,4], \,[5, \,6], \,[5, 
\,7], \,[7, \,2], \,[7, \,4], \,[7, \,6]] }
}
\end{maplelatex}

\end{maplegroup}
\begin{maplegroup}
We list here all intervals of points in $P$ which are admitted by the
given ordering $>$:

\end{maplegroup}
\begin{maplegroup}
\begin{_pstyle1}
\begin{mapleinput}
\mapleinline{active}{1d}{I1:=GenerateIntervals(P,rel);}{%
}
\end{mapleinput}
\end{_pstyle1}

\mapleresult
\begin{maplelatex}
\mapleinline{inert}{2d}{I1 := [[], [2], [1], [0], [7], [3], [4], [5], [6], [3, 5], [2, 6],
[0, 7], [0, 2], [1, 7], [3, 7], [5, 7], [0, 4], [0, 6], [0, 1], [0,
3], [0, 5], [1, 3], [1, 5], [2, 4], [2, 7], [4, 6], [4, 7], [6, 7],
[1, 3, 5], [2, 4, 6], [0, 6, 7], [2, 6, 7], [4, 6, 7], [0, 1, 7], [0,
2, 7], [0, 3, 7], [1, 3, 7], [0, 4, 7], [2, 4, 7], [0, 5, 7], [1, 5,
7], [3, 5, 7], [0, 1, 3], [0, 2, 4], [0, 1, 5], [0, 3, 5], [0, 2, 6],
[0, 4, 6], [0, 5, 6, 7], [0, 1, 6, 7], [0, 2, 6, 7], [0, 3, 6, 7], [0,
4, 6, 7], [2, 4, 6, 7], [0, 1, 2, 7], [0, 1, 3, 7], [0, 2, 3, 7], [0,
1, 4, 7], [0, 2, 4, 7], [0, 3, 4, 7], [0, 1, 5, 7], [0, 2, 5, 7], [0,
3, 5, 7], [1, 3, 5, 7], [0, 4, 5, 7], [0, 1, 3, 5], [0, 2, 4, 6], [0,
4, 5, 6, 7], [0, 1, 5, 6, 7], [0, 2, 5, 6, 7], [0, 3, 5, 6, 7], [0, 1,
2, 6, 7], [0, 1, 3, 6, 7], [0, 2, 3, 6, 7], [0, 1, 4, 6, 7], [0, 2, 4,
6, 7], [0, 3, 4, 6, 7], [0, 1, 2, 3, 7], [0, 1, 2, 4, 7], [0, 1, 3, 4,
7], [0, 2, 3, 4, 7], [0, 1, 2, 5, 7], [0, 1, 3, 5, 7], [0, 2, 3, 5,
7], [0, 1, 4, 5, 7], [0, 2, 4, 5, 7], [0, 3, 4, 5, 7], [0, 3, 4, 5, 6,
7], [0, 1, 4, 5, 6, 7], [0, 2, 4, 5, 6, 7], [0, 1, 2, 5, 6, 7], [0, 1,
3, 5, 6, 7], [0, 2, 3, 5, 6, 7], [0, 1, 2, 3, 6, 7], [0, 1, 2, 4, 6,
7], [0, 1, 3, 4, 6, 7], [0, 2, 3, 4, 6, 7], [0, 1, 2, 3, 4, 7], [0, 1,
2, 3, 5, 7], [0, 1, 2, 4, 5, 7], [0, 1, 3, 4, 5, 7], [0, 2, 3, 4, 5,
7], [0, 2, 3, 4, 5, 6, 7], [0, 1, 3, 4, 5, 6, 7], [0, 1, 2, 4, 5, 6,
7], [0, 1, 2, 3, 5, 6, 7], [0, 1, 2, 3, 4, 6, 7], [0, 1, 2, 3, 4, 5,
7], [0, 1, 2, 3, 4, 5, 6, 7]];}{%
\maplemultiline{
\mathit{I1} := [[], \,[2], \,[1], \,[0], \,[7], \,[3], \,[4], \,[
5], \,[6], \,[3, \,5], \,[2, \,6], \,[0, \,7], \,[0, \,2], \,[1, 
\,7], \,[3, \,7], \,[5, \,7],  \\
[0, \,4], \,[0, \,6], \,[0, \,1], \,[0, \,3], \,[0, \,5], \,[1, 
\,3], \,[1, \,5], \,[2, \,4], \,[2, \,7], \,[4, \,6], \,[4, \,7]
, \,[6, \,7],  \\
[1, \,3, \,5], \,[2, \,4, \,6], \,[0, \,6, \,7], \,[2, \,6, \,7]
, \,[4, \,6, \,7], \,[0, \,1, \,7], \,[0, \,2, \,7], \,[0, \,3, 
\,7], \,[1, \,3, \,7],  \\
[0, \,4, \,7], \,[2, \,4, \,7], \,[0, \,5, \,7], \,[1, \,5, \,7]
, \,[3, \,5, \,7], \,[0, \,1, \,3], \,[0, \,2, \,4], \,[0, \,1, 
\,5], \,[0, \,3, \,5],  \\
[0, \,2, \,6], \,[0, \,4, \,6], \,[0, \,5, \,6, \,7], \,[0, \,1, 
\,6, \,7], \,[0, \,2, \,6, \,7], \,[0, \,3, \,6, \,7], \,[0, \,4
, \,6, \,7],  \\
[2, \,4, \,6, \,7], \,[0, \,1, \,2, \,7], \,[0, \,1, \,3, \,7], 
\,[0, \,2, \,3, \,7], \,[0, \,1, \,4, \,7], \,[0, \,2, \,4, \,7]
, \,[0, \,3, \,4, \,7],  \\
[0, \,1, \,5, \,7], \,[0, \,2, \,5, \,7], \,[0, \,3, \,5, \,7], 
\,[1, \,3, \,5, \,7], \,[0, \,4, \,5, \,7], \,[0, \,1, \,3, \,5]
, \,[0, \,2, \,4, \,6],  \\
[0, \,4, \,5, \,6, \,7], \,[0, \,1, \,5, \,6, \,7], \,[0, \,2, \,
5, \,6, \,7], \,[0, \,3, \,5, \,6, \,7], \,[0, \,1, \,2, \,6, \,7
],  \\
[0, \,1, \,3, \,6, \,7], \,[0, \,2, \,3, \,6, \,7], \,[0, \,1, \,
4, \,6, \,7], \,[0, \,2, \,4, \,6, \,7], \,[0, \,3, \,4, \,6, \,7
],  \\
[0, \,1, \,2, \,3, \,7], \,[0, \,1, \,2, \,4, \,7], \,[0, \,1, \,
3, \,4, \,7], \,[0, \,2, \,3, \,4, \,7], \,[0, \,1, \,2, \,5, \,7
],  \\
[0, \,1, \,3, \,5, \,7], \,[0, \,2, \,3, \,5, \,7], \,[0, \,1, \,
4, \,5, \,7], \,[0, \,2, \,4, \,5, \,7], \,[0, \,3, \,4, \,5, \,7
],  \\
[0, \,3, \,4, \,5, \,6, \,7], \,[0, \,1, \,4, \,5, \,6, \,7], \,[
0, \,2, \,4, \,5, \,6, \,7], \,[0, \,1, \,2, \,5, \,6, \,7], \,[0
, \,1, \,3, \,5, \,6, \,7],  \\
[0, \,2, \,3, \,5, \,6, \,7], \,[0, \,1, \,2, \,3, \,6, \,7], \,[
0, \,1, \,2, \,4, \,6, \,7], \,[0, \,1, \,3, \,4, \,6, \,7], \,[0
, \,2, \,3, \,4, \,6, \,7],  \\
[0, \,1, \,2, \,3, \,4, \,7], \,[0, \,1, \,2, \,3, \,5, \,7], \,[
0, \,1, \,2, \,4, \,5, \,7], \,[0, \,1, \,3, \,4, \,5, \,7], \,[0
, \,2, \,3, \,4, \,5, \,7],  \\
[0, \,2, \,3, \,4, \,5, \,6, \,7], \,[0, \,1, \,3, \,4, \,5, \,6
, \,7], \,[0, \,1, \,2, \,4, \,5, \,6, \,7], \,[0, \,1, \,2, \,3
, \,5, \,6, \,7],  \\
[0, \,1, \,2, \,3, \,4, \,6, \,7], \,[0, \,1, \,2, \,3, \,4, \,5
, \,7], \,[0, \,1, \,2, \,3, \,4, \,5, \,6, \,7]] }
}
\end{maplelatex}

\end{maplegroup}
\begin{maplegroup}
\begin{_pstyle1}
\begin{mapleinput}
\mapleinline{active}{1d}{#I123:=GenerateAdjacentIntervals(P,rel);}{%
}
\end{mapleinput}
\end{_pstyle1}

\end{maplegroup}
\begin{maplegroup}
\begin{_pstyle1}
\begin{mapleinput}
\mapleinline{active}{1d}{orbits:=Orbits(D6C2,ord);}{%
}
\end{mapleinput}
\end{_pstyle1}

\mapleresult
\begin{maplelatex}
\mapleinline{inert}{2d}{orbits := [[[0, 2], [0, 4], [0, 6], [7, 2], [7, 4], [7, 6]], [[1, 0],
[3, 0], [1, 7], [5, 0], [3, 7], [5, 7]], [[1, 2], [3, 4], [1, 6], [5,
6], [5, 4], [3, 2]], [[1, 4], [3, 6], [5, 2]]];}{%
\maplemultiline{
\mathit{orbits} := [[[0, \,2], \,[0, \,4], \,[0, \,6], \,[7, \,2]
, \,[7, \,4], \,[7, \,6]],  \\
[[1, \,0], \,[3, \,0], \,[1, \,7], \,[5, \,0], \,[3, \,7], \,[5, 
\,7]], \,[[1, \,2], \,[3, \,4], \,[1, \,6], \,[5, \,6], \,[5, \,4
], \,[3, \,2]],  \\
[[1, \,4], \,[3, \,6], \,[5, \,2]]] }
}
\end{maplelatex}

\end{maplegroup}
\begin{maplegroup}
Next we provide the homology {\sc Conley} indices. Note that since $0$
and $7$ are non-hyperbolic, their homology is not of rank $1$.

\end{maplegroup}
\begin{maplegroup}
\begin{_pstyle1}
\begin{mapleinput}
\mapleinline{active}{1d}{CHp:=[[0,2],2,0,2,0,2,0,[0,2]];}{%
}
\end{mapleinput}
\end{_pstyle1}

\mapleresult
\begin{maplelatex}
\mapleinline{inert}{2d}{CHp := [[0, 2], 2, 0, 2, 0, 2, 0, [0, 2]];}{%
\[
\mathit{CHp} := [[0, \,2], \,2, \,0, \,2, \,0, \,2, \,0, \,[0, \,
2]]
\]
}
\end{maplelatex}

\end{maplegroup}
\begin{maplegroup}
\begin{mapleinput}
\mapleinline{active}{1d}{CH:=[op(zip((x,y)->x=y,P,CHp)),P=[1,0,1]];}{%
}
\end{mapleinput}

\mapleresult
\begin{maplelatex}
\mapleinline{inert}{2d}{CH := [0 = [0, 2], 1 = 2, 2 = 0, 3 = 2, 4 = 0, 5 = 2, 6 = 0, 7 = [0,
2], [0, 1, 2, 3, 4, 5, 6, 7] = [1, 0, 1]];}{%
\maplemultiline{
\mathit{CH} := [0=[0, \,2], \,1=2, \,2=0, \,3=2, \,4=0, \,5=2, \,
6=0, \,7=[0, \,2],  \\
[0, \,1, \,2, \,3, \,4, \,5, \,6, \,7]=[1, \,0, \,1]] }
}
\end{maplelatex}

\end{maplegroup}
\begin{maplegroup}
Now we compute all connection matrices for the given homology {\sc
Conley}
indices respecting the symmetry group $\Gamma$ generated by
$\mathtt{D6C2}$,
where the action of $\Gamma$ on $\bigoplus_{p \in P} CH_i(p)$ for $i =
0, 1,
2$ by the matrices defined above is taken into account. The optional
string
$\mathtt{"std\_notation"}$ means that the matrix
$\Delta_i:\bigoplus_{p \in P} CH_i(p) \to \bigoplus_{p \in P}
CH_{i-1}(p)$ is
a $d_i\times d_{i-1}$-matrix, where $d_j:=\dim_{\F_2}\bigoplus_{p \in
P} CH_j(p)$:

\end{maplegroup}
\begin{maplegroup}
\begin{_pstyle1}
\begin{mapleinput}
\mapleinline{active}{1d}{Con:=ConnectionMatrices(P,rel,CH,var,"Symmetry"=D6C2_rep,"std_notatio
n");}{%
}
\end{mapleinput}
\end{_pstyle1}

\mapleresult
\begin{maplelatex}
\mapleinline{inert}{2d}{Con := [[matrix([[0, 1, 1], [1, 0, 1], [0, 1, 1], [1, 0, 1]]),
matrix([[0, 1, 0, 1], [1, 1, 1, 1], [1, 0, 1, 0]])]];}{%
\[
\mathit{Con} :=  \left[  \!  \left[  \!  \left[ 
{\begin{array}{rrr}
0 & 1 & 1 \\
1 & 0 & 1 \\
0 & 1 & 1 \\
1 & 0 & 1
\end{array}}
 \right] , \, \left[ 
{\begin{array}{rrrr}
0 & 1 & 0 & 1 \\
1 & 1 & 1 & 1 \\
1 & 0 & 1 & 0
\end{array}}
 \right]  \!  \right]  \!  \right] 
\]
}
\end{maplelatex}

\end{maplegroup}
\begin{maplegroup}
We conclude that there is exactly one connection matrix which is
compatible with the given data. It is given by the list containing the
matrices which represent the boundary maps from $\bigoplus_{p \in P}
CH_1(p)$ to $\bigoplus_{p \in P} CH_0(p)$ resp.\  from $\bigoplus_{p
\in P} CH_2(p)$ to $\bigoplus_{p \in P} CH_1(p)$. Here we extract the
information about connecting orbits between points in $P$ which is
provided by the computed connection matrix:

\end{maplegroup}
\begin{maplegroup}
\begin{mapleinput}
\mapleinline{active}{1d}{ConnectionSubmatrix(Con[1],P,[0,2],CHp,var),
ConnectionSubmatrix(Con[1],P,[0,4],CHp,var),
ConnectionSubmatrix(Con[1],P,[0,6],CHp,var);
ConnectionSubmatrix(Con[1],P,[7,2],CHp,var),
ConnectionSubmatrix(Con[1],P,[7,4],CHp,var),
ConnectionSubmatrix(Con[1],P,[7,6],CHp,var);}{%
}
\end{mapleinput}

\mapleresult
\begin{maplelatex}
\mapleinline{inert}{2d}{[matrix([[0], [1]]), matrix([[0]])], [matrix([[1], [0]]),
matrix([[0]])], [matrix([[1], [1]]), matrix([[0]])];}{%
\[
[ \left[ 
{\begin{array}{r}
0 \\
1
\end{array}}
 \right] , \, \left[ 
{\begin{array}{r}
0
\end{array}}
 \right] ], \,[ \left[ 
{\begin{array}{r}
1 \\
0
\end{array}}
 \right] , \, \left[ 
{\begin{array}{r}
0
\end{array}}
 \right] ], \,[ \left[ 
{\begin{array}{r}
1 \\
1
\end{array}}
 \right] , \, \left[ 
{\begin{array}{r}
0
\end{array}}
 \right] ]
\]
}
\end{maplelatex}

\begin{maplelatex}
\mapleinline{inert}{2d}{[matrix([[0], [1]]), matrix([[0]])], [matrix([[1], [0]]),
matrix([[0]])], [matrix([[1], [1]]), matrix([[0]])];}{%
\[
[ \left[ 
{\begin{array}{r}
0 \\
1
\end{array}}
 \right] , \, \left[ 
{\begin{array}{r}
0
\end{array}}
 \right] ], \,[ \left[ 
{\begin{array}{r}
1 \\
0
\end{array}}
 \right] , \, \left[ 
{\begin{array}{r}
0
\end{array}}
 \right] ], \,[ \left[ 
{\begin{array}{r}
1 \\
1
\end{array}}
 \right] , \, \left[ 
{\begin{array}{r}
0
\end{array}}
 \right] ]
\]
}
\end{maplelatex}

\end{maplegroup}
\begin{maplegroup}
Hence, there are connecting orbits from $0$ to $2$, to $4$, to $6$,
and from $7$ to $2$, to $4$, and to $6$.

\end{maplegroup}
\begin{maplegroup}
\begin{mapleinput}
\mapleinline{active}{1d}{ConnectionSubmatrix(Con[1],P,[1,0],CHp,var),
ConnectionSubmatrix(Con[1],P,[1,7],CHp,var);
ConnectionSubmatrix(Con[1],P,[3,0],CHp,var),
ConnectionSubmatrix(Con[1],P,[3,7],CHp,var);
ConnectionSubmatrix(Con[1],P,[5,0],CHp,var),
ConnectionSubmatrix(Con[1],P,[5,7],CHp,var);}{%
}
\end{mapleinput}

\mapleresult
\begin{maplelatex}
\mapleinline{inert}{2d}{[matrix([[0]]), matrix([[0, 1]])], [matrix([[0]]), matrix([[0,
1]])];}{%
\[
[ \left[ 
{\begin{array}{r}
0
\end{array}}
 \right] , \, \left[ 
{\begin{array}{rr}
0 & 1
\end{array}}
 \right] ], \,[ \left[ 
{\begin{array}{r}
0
\end{array}}
 \right] , \, \left[ 
{\begin{array}{rr}
0 & 1
\end{array}}
 \right] ]
\]
}
\end{maplelatex}

\begin{maplelatex}
\mapleinline{inert}{2d}{[matrix([[0]]), matrix([[1, 1]])], [matrix([[0]]), matrix([[1,
1]])];}{%
\[
[ \left[ 
{\begin{array}{r}
0
\end{array}}
 \right] , \, \left[ 
{\begin{array}{rr}
1 & 1
\end{array}}
 \right] ], \,[ \left[ 
{\begin{array}{r}
0
\end{array}}
 \right] , \, \left[ 
{\begin{array}{rr}
1 & 1
\end{array}}
 \right] ]
\]
}
\end{maplelatex}

\begin{maplelatex}
\mapleinline{inert}{2d}{[matrix([[0]]), matrix([[1, 0]])], [matrix([[0]]), matrix([[1,
0]])];}{%
\[
[ \left[ 
{\begin{array}{r}
0
\end{array}}
 \right] , \, \left[ 
{\begin{array}{rr}
1 & 0
\end{array}}
 \right] ], \,[ \left[ 
{\begin{array}{r}
0
\end{array}}
 \right] , \, \left[ 
{\begin{array}{rr}
1 & 0
\end{array}}
 \right] ]
\]
}
\end{maplelatex}

\end{maplegroup}
\begin{maplegroup}
Similarly, there are connecting orbits from $1$ to $0$ and $7$, from
$3$ to $0$ and $7$, and from $5$ to $0$ and $7$.

\end{maplegroup}
%
\POP

\subsection{The {\sc Cahn-Hilliard} equation}

\PUSH{Cahn-Hilliard.tex}%
\DefineParaStyle{Maple Output}
\DefineCharStyle{2D Math}
\DefineCharStyle{2D Output}
\DefineCharStyle{LaTeX}
\begin{maplegroup}
In this example we consider a $D_8$-symmetric flow with nine
hyperbolic equilibria studied in {\cite[{2.2, 4.2}]{MMW}}. We stick to
the notation of the equilibria of that reference: 
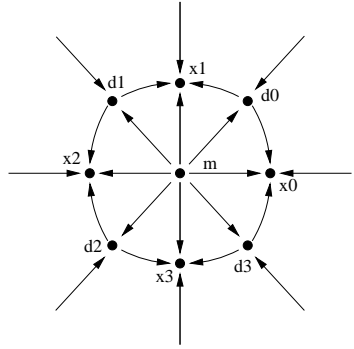
\begin{center}
\includegraphics[scale=0.6]{Cahn-Hilliard.pstex} \end{center}

\end{maplegroup}
\begin{maplegroup}
\begin{mapleinput}
\mapleinline{active}{1d}{restart;}{%
}
\end{mapleinput}

\end{maplegroup}
\begin{maplegroup}
\begin{mapleinput}
\mapleinline{active}{1d}{with(conley): with(PIR): with(homalg):}{%
}
\end{mapleinput}

\end{maplegroup}
\begin{maplegroup}
\begin{mapleinput}
\mapleinline{active}{1d}{`homalg/default`:=`PIR/homalg`:}{%
}
\end{mapleinput}

\end{maplegroup}
\begin{maplegroup}
\begin{mapleinput}
\mapleinline{active}{1d}{var:=[2];}{%
}
\end{mapleinput}

\mapleresult
\begin{maplelatex}
\mapleinline{inert}{2d}{var := [2];}{%
\[
\mathit{var} := [2]
\]
}
\end{maplelatex}

\end{maplegroup}
\begin{maplegroup}
\begin{mapleinput}
\mapleinline{active}{1d}{Pvar(var);}{%
}
\end{mapleinput}

\mapleresult
\begin{maplelatex}
\mapleinline{inert}{2d}{["Z/pZ", 2];}{%
\[
[\mbox{``Z/pZ''}, \,2]
\]
}
\end{maplelatex}

\end{maplegroup}
\begin{maplegroup}
\begin{mapleinput}
\mapleinline{active}{1d}{P:=[x0,x1,x2,x3,d0,d1,d2,d3,m];}{%
}
\end{mapleinput}

\mapleresult
\begin{maplelatex}
\mapleinline{inert}{2d}{P := [x0, x1, x2, x3, d0, d1, d2, d3, m];}{%
\[
P := [\mathit{x0}, \,\mathit{x1}, \,\mathit{x2}, \,\mathit{x3}, 
\,\mathit{d0}, \,\mathit{d1}, \,\mathit{d2}, \,\mathit{d3}, \,m]
\]
}
\end{maplelatex}

\end{maplegroup}
\begin{maplegroup}
\begin{mapleinput}
\mapleinline{active}{1d}{CHp:=[0,0,0,0,1,1,1,1,2];}{%
}
\end{mapleinput}

\mapleresult
\begin{maplelatex}
\mapleinline{inert}{2d}{CHp := [0, 0, 0, 0, 1, 1, 1, 1, 2];}{%
\[
\mathit{CHp} := [0, \,0, \,0, \,0, \,1, \,1, \,1, \,1, \,2]
\]
}
\end{maplelatex}

\end{maplegroup}
\begin{maplegroup}
\begin{mapleinput}
\mapleinline{active}{1d}{rel:=DefinePartialOrderByPotential(P,CHp,P);}{%
}
\end{mapleinput}

\mapleresult
\begin{maplelatex}
\mapleinline{inert}{2d}{rel := [[d0, x0], [d1, x0], [d2, x0], [d3, x0], [m, x0], [d0, x1],
[d1, x1], [d2, x1], [d3, x1], [m, x1], [d0, x2], [d1, x2], [d2, x2],
[d3, x2], [m, x2], [d0, x3], [d1, x3], [d2, x3], [d3, x3], [m, x3],
[m, d0], [m, d1], [m, d2], [m, d3]];}{%
\maplemultiline{
\mathit{rel} := [[\mathit{d0}, \,\mathit{x0}], \,[\mathit{d1}, \,
\mathit{x0}], \,[\mathit{d2}, \,\mathit{x0}], \,[\mathit{d3}, \,
\mathit{x0}], \,[m, \,\mathit{x0}], \,[\mathit{d0}, \,\mathit{x1}
], \,[\mathit{d1}, \,\mathit{x1}], \,[\mathit{d2}, \,\mathit{x1}]
, \,[\mathit{d3}, \,\mathit{x1}],  \\
[m, \,\mathit{x1}], \,[\mathit{d0}, \,\mathit{x2}], \,[\mathit{d1
}, \,\mathit{x2}], \,[\mathit{d2}, \,\mathit{x2}], \,[\mathit{d3}
, \,\mathit{x2}], \,[m, \,\mathit{x2}], \,[\mathit{d0}, \,
\mathit{x3}], \,[\mathit{d1}, \,\mathit{x3}], \,[\mathit{d2}, \,
\mathit{x3}],  \\
[\mathit{d3}, \,\mathit{x3}], \,[m, \,\mathit{x3}], \,[m, \,
\mathit{d0}], \,[m, \,\mathit{d1}], \,[m, \,\mathit{d2}], \,[m, 
\,\mathit{d3}]] }
}
\end{maplelatex}

\end{maplegroup}
\begin{maplegroup}
\begin{mapleinput}
\mapleinline{active}{1d}{ord:=GeneratePartialOrder(P,rel);}{%
}
\end{mapleinput}

\mapleresult
\begin{maplelatex}
\mapleinline{inert}{2d}{ord := [[d0, x0], [d0, x1], [d0, x2], [d0, x3], [d1, x0], [d1, x1],
[d1, x2], [d1, x3], [d2, x0], [d2, x1], [d2, x2], [d2, x3], [d3, x0],
[d3, x1], [d3, x2], [d3, x3], [m, x0], [m, x1], [m, x2], [m, x3], [m,
d0], [m, d1], [m, d2], [m, d3]];}{%
\maplemultiline{
\mathit{ord} := [[\mathit{d0}, \,\mathit{x0}], \,[\mathit{d0}, \,
\mathit{x1}], \,[\mathit{d0}, \,\mathit{x2}], \,[\mathit{d0}, \,
\mathit{x3}], \,[\mathit{d1}, \,\mathit{x0}], \,[\mathit{d1}, \,
\mathit{x1}], \,[\mathit{d1}, \,\mathit{x2}], \,[\mathit{d1}, \,
\mathit{x3}], \,[\mathit{d2}, \,\mathit{x0}],  \\
[\mathit{d2}, \,\mathit{x1}], \,[\mathit{d2}, \,\mathit{x2}], \,[
\mathit{d2}, \,\mathit{x3}], \,[\mathit{d3}, \,\mathit{x0}], \,[
\mathit{d3}, \,\mathit{x1}], \,[\mathit{d3}, \,\mathit{x2}], \,[
\mathit{d3}, \,\mathit{x3}], \,[m, \,\mathit{x0}], \,[m, \,
\mathit{x1}],  \\
[m, \,\mathit{x2}], \,[m, \,\mathit{x3}], \,[m, \,\mathit{d0}], 
\,[m, \,\mathit{d1}], \,[m, \,\mathit{d2}], \,[m, \,\mathit{d3}]]
 }
}
\end{maplelatex}

\end{maplegroup}
\begin{maplegroup}
\begin{mapleinput}
\mapleinline{active}{1d}{D8:=[ [[x0,x1,x2,x3],[d0,d1,d2,d3]], [[x1,x3],[d0,d3],[d1,d2]] ];}{%
}
\end{mapleinput}

\mapleresult
\begin{maplelatex}
\mapleinline{inert}{2d}{D8 := [[[x0, x1, x2, x3], [d0, d1, d2, d3]], [[x1, x3], [d0, d3],
[d1, d2]]];}{%
\[
\mathrm{D8} := [[[\mathit{x0}, \,\mathit{x1}, \,\mathit{x2}, \,
\mathit{x3}], \,[\mathit{d0}, \,\mathit{d1}, \,\mathit{d2}, \,
\mathit{d3}]], \,[[\mathit{x1}, \,\mathit{x3}], \,[\mathit{d0}, 
\,\mathit{d3}], \,[\mathit{d1}, \,\mathit{d2}]]]
\]
}
\end{maplelatex}

\end{maplegroup}
\begin{maplegroup}
\begin{mapleinput}
\mapleinline{active}{1d}{CH:=[op(zip((x,y)->x=y,P,CHp)),P=0];}{%
}
\end{mapleinput}

\mapleresult
\begin{maplelatex}
\mapleinline{inert}{2d}{CH := [x0 = 0, x1 = 0, x2 = 0, x3 = 0, d0 = 1, d1 = 1, d2 = 1, d3 =
1, m = 2, [x0, x1, x2, x3, d0, d1, d2, d3, m] = 0];}{%
\maplemultiline{
\mathit{CH} := [\mathit{x0}=0, \,\mathit{x1}=0, \,\mathit{x2}=0, 
\,\mathit{x3}=0, \,\mathit{d0}=1, \,\mathit{d1}=1, \,\mathit{d2}=
1, \,\mathit{d3}=1, \,m=2,  \\
[\mathit{x0}, \,\mathit{x1}, \,\mathit{x2}, \,\mathit{x3}, \,
\mathit{d0}, \,\mathit{d1}, \,\mathit{d2}, \,\mathit{d3}, \,m]=0]
 }
}
\end{maplelatex}

\end{maplegroup}
\begin{maplegroup}
\begin{mapleinput}
\mapleinline{active}{1d}{Con:=ConnectionMatrices(P,rel,CH,var,"Symmetry"=D8,"full_notation");}
{%
}
\end{mapleinput}

\mapleresult
\begin{maplelatex}
\mapleinline{inert}{2d}{Con := [[matrix([[0, 0, 0, 0, 0, 0, 0, 0, 0], [0, 0, 0, 0, 0, 0, 0,
0, 0], [0, 0, 0, 0, 0, 0, 0, 0, 0], [0, 0, 0, 0, 0, 0, 0, 0, 0], [1,
1, 0, 0, 0, 0, 0, 0, 0], [0, 1, 1, 0, 0, 0, 0, 0, 0], [0, 0, 1, 1, 0,
0, 0, 0, 0], [1, 0, 0, 1, 0, 0, 0, 0, 0], [0, 0, 0, 0, 0, 0, 0, 0,
0]]), matrix([[0, 0, 0, 0, 0, 0, 0, 0, 0], [0, 0, 0, 0, 0, 0, 0, 0,
0], [0, 0, 0, 0, 0, 0, 0, 0, 0], [0, 0, 0, 0, 0, 0, 0, 0, 0], [0, 0,
0, 0, 0, 0, 0, 0, 0], [0, 0, 0, 0, 0, 0, 0, 0, 0], [0, 0, 0, 0, 0, 0,
0, 0, 0], [0, 0, 0, 0, 0, 0, 0, 0, 0], [0, 0, 0, 0, 1, 1, 1, 1, 0]])],
[matrix([[0, 0, 0, 0, 0, 0, 0, 0, 0], [0, 0, 0, 0, 0, 0, 0, 0, 0], [0,
0, 0, 0, 0, 0, 0, 0, 0], [0, 0, 0, 0, 0, 0, 0, 0, 0], [0, 0, 1, 1, 0,
0, 0, 0, 0], [1, 0, 0, 1, 0, 0, 0, 0, 0], [1, 1, 0, 0, 0, 0, 0, 0, 0],
[0, 1, 1, 0, 0, 0, 0, 0, 0], [0, 0, 0, 0, 0, 0, 0, 0, 0]]),
matrix([[0, 0, 0, 0, 0, 0, 0, 0, 0], [0, 0, 0, 0, 0, 0, 0, 0, 0], [0,
0, 0, 0, 0, 0, 0, 0, 0], [0, 0, 0, 0, 0, 0, 0, 0, 0], [0, 0, 0, 0, 0,
0, 0, 0, 0], [0, 0, 0, 0, 0, 0, 0, 0, 0], [0, 0, 0, 0, 0, 0, 0, 0, 0],
[0, 0, 0, 0, 0, 0, 0, 0, 0], [0, 0, 0, 0, 1, 1, 1, 1, 0]])]];}{%
\maplemultiline{
\mathit{Con} :=  \left[ {\vrule height3.31em width0em depth3.31em
} \right. \!  \!  \left[  \!  \left[ 
{\begin{array}{rrrrrrrrr}
0 & 0 & 0 & 0 & 0 & 0 & 0 & 0 & 0 \\
0 & 0 & 0 & 0 & 0 & 0 & 0 & 0 & 0 \\
0 & 0 & 0 & 0 & 0 & 0 & 0 & 0 & 0 \\
0 & 0 & 0 & 0 & 0 & 0 & 0 & 0 & 0 \\
1 & 1 & 0 & 0 & 0 & 0 & 0 & 0 & 0 \\
0 & 1 & 1 & 0 & 0 & 0 & 0 & 0 & 0 \\
0 & 0 & 1 & 1 & 0 & 0 & 0 & 0 & 0 \\
1 & 0 & 0 & 1 & 0 & 0 & 0 & 0 & 0 \\
0 & 0 & 0 & 0 & 0 & 0 & 0 & 0 & 0
\end{array}}
 \right] , \, \left[ 
{\begin{array}{rrrrrrrrr}
0 & 0 & 0 & 0 & 0 & 0 & 0 & 0 & 0 \\
0 & 0 & 0 & 0 & 0 & 0 & 0 & 0 & 0 \\
0 & 0 & 0 & 0 & 0 & 0 & 0 & 0 & 0 \\
0 & 0 & 0 & 0 & 0 & 0 & 0 & 0 & 0 \\
0 & 0 & 0 & 0 & 0 & 0 & 0 & 0 & 0 \\
0 & 0 & 0 & 0 & 0 & 0 & 0 & 0 & 0 \\
0 & 0 & 0 & 0 & 0 & 0 & 0 & 0 & 0 \\
0 & 0 & 0 & 0 & 0 & 0 & 0 & 0 & 0 \\
0 & 0 & 0 & 0 & 1 & 1 & 1 & 1 & 0
\end{array}}
 \right]  \!  \right] ,  \\
 \left[  \!  \left[ 
{\begin{array}{rrrrrrrrr}
0 & 0 & 0 & 0 & 0 & 0 & 0 & 0 & 0 \\
0 & 0 & 0 & 0 & 0 & 0 & 0 & 0 & 0 \\
0 & 0 & 0 & 0 & 0 & 0 & 0 & 0 & 0 \\
0 & 0 & 0 & 0 & 0 & 0 & 0 & 0 & 0 \\
0 & 0 & 1 & 1 & 0 & 0 & 0 & 0 & 0 \\
1 & 0 & 0 & 1 & 0 & 0 & 0 & 0 & 0 \\
1 & 1 & 0 & 0 & 0 & 0 & 0 & 0 & 0 \\
0 & 1 & 1 & 0 & 0 & 0 & 0 & 0 & 0 \\
0 & 0 & 0 & 0 & 0 & 0 & 0 & 0 & 0
\end{array}}
 \right] , \, \left[ 
{\begin{array}{rrrrrrrrr}
0 & 0 & 0 & 0 & 0 & 0 & 0 & 0 & 0 \\
0 & 0 & 0 & 0 & 0 & 0 & 0 & 0 & 0 \\
0 & 0 & 0 & 0 & 0 & 0 & 0 & 0 & 0 \\
0 & 0 & 0 & 0 & 0 & 0 & 0 & 0 & 0 \\
0 & 0 & 0 & 0 & 0 & 0 & 0 & 0 & 0 \\
0 & 0 & 0 & 0 & 0 & 0 & 0 & 0 & 0 \\
0 & 0 & 0 & 0 & 0 & 0 & 0 & 0 & 0 \\
0 & 0 & 0 & 0 & 0 & 0 & 0 & 0 & 0 \\
0 & 0 & 0 & 0 & 1 & 1 & 1 & 1 & 0
\end{array}}
 \right]  \!  \right]  \! \! \left. {\vrule 
height3.31em width0em depth3.31em} \right]  }
}
\end{maplelatex}

\end{maplegroup}
\begin{maplegroup}
We obtain two connection matrices. The full flow diagram above would
give a unique $D_8$-symmetric connection matrix, but this diagram was
derived in \cite{MMW} using dynamical arguments. Here we demonstrated
that the lack of the dynamical arguments produces ambiguities. This is
the point of departure in a joint work in progress with {\sc
Maier-Paape}.

\end{maplegroup}
%
\POP

\subsection{An example for a $C$-connection matrix}

\PUSH{CCon.tex}%
\DefineParaStyle{Maple Output}
\DefineCharStyle{2D Math}
\DefineCharStyle{2D Output}
\DefineCharStyle{LaTeX}
\begin{maplegroup}
Here we demonstrate on a tiny example the advantage of a
$C$-connection matrix over a connection matrix. The aim is to look at
the two equilibria $0$ and $1$ of the flow visualized in
\begin{center}  \includegraphics[scale=0.6]{1dim1.pstex} \end{center}
as a black box. Then the flow can also be depicted as follows:
\begin{center}  \includegraphics[scale=0.6]{1dim2.pstex} \end{center}
Here we named the black box equilibrium again by $1$.

\medskip

\end{maplegroup}
\begin{maplegroup}
\begin{mapleinput}
\mapleinline{active}{1d}{restart;}{%
}
\end{mapleinput}

\end{maplegroup}
\begin{maplegroup}
\begin{mapleinput}
\mapleinline{active}{1d}{with(conley): with(PIR): with(homalg):}{%
}
\end{mapleinput}

\end{maplegroup}
\begin{maplegroup}
\begin{mapleinput}
\mapleinline{active}{1d}{`homalg/default`:=`PIR/homalg`:}{%
}
\end{mapleinput}

\end{maplegroup}
\begin{maplegroup}
We set up the field of coefficients to be $\F_2:=\Z/2\Z$:

\end{maplegroup}
\begin{maplegroup}
\begin{mapleinput}
\mapleinline{active}{1d}{var:=[2];}{%
}
\end{mapleinput}

\mapleresult
\begin{maplelatex}
\mapleinline{inert}{2d}{var := [2];}{%
\[
\mathit{var} := [2]
\]
}
\end{maplelatex}

\end{maplegroup}
\begin{maplegroup}
\begin{mapleinput}
\mapleinline{active}{1d}{Pvar(var);}{%
}
\end{mapleinput}

\mapleresult
\begin{maplelatex}
\mapleinline{inert}{2d}{["Z/pZ", 2];}{%
\[
[\mbox{``Z/pZ''}, \,2]
\]
}
\end{maplelatex}

\end{maplegroup}
\begin{maplegroup}
\begin{mapleinput}
\mapleinline{active}{1d}{P:=[0,1,2];}{%
}
\end{mapleinput}

\mapleresult
\begin{maplelatex}
\mapleinline{inert}{2d}{P := [0, 1, 2];}{%
\[
P := [0, \,1, \,2]
\]
}
\end{maplelatex}

\end{maplegroup}
\begin{maplegroup}
The potential values/or the indices:

\end{maplegroup}
\begin{maplegroup}
\begin{mapleinput}
\mapleinline{active}{1d}{CHp:=[1,0,1];}{%
}
\end{mapleinput}

\mapleresult
\begin{maplelatex}
\mapleinline{inert}{2d}{CHp := [1, 0, 1];}{%
\[
\mathit{CHp} := [1, \,0, \,1]
\]
}
\end{maplelatex}

\end{maplegroup}
\begin{maplegroup}
\begin{mapleinput}
\mapleinline{active}{1d}{rel:=[[0,1],[2,1]];}{%
}
\end{mapleinput}

\mapleresult
\begin{maplelatex}
\mapleinline{inert}{2d}{rel := [[0, 1], [2, 1]];}{%
\[
\mathit{rel} := [[0, \,1], \,[2, \,1]]
\]
}
\end{maplelatex}

\end{maplegroup}
\begin{maplegroup}
\begin{mapleinput}
\mapleinline{active}{1d}{ord:=GeneratePartialOrder(P,rel);}{%
}
\end{mapleinput}

\mapleresult
\begin{maplelatex}
\mapleinline{inert}{2d}{ord := [[0, 1], [2, 1]];}{%
\[
\mathit{ord} := [[0, \,1], \,[2, \,1]]
\]
}
\end{maplelatex}

\end{maplegroup}
\begin{maplegroup}
\begin{mapleinput}
\mapleinline{active}{1d}{I1:=GenerateIntervals(P,rel);}{%
}
\end{mapleinput}

\mapleresult
\begin{maplelatex}
\mapleinline{inert}{2d}{I1 := [[], [2], [0], [1], [0, 1], [0, 2], [1, 2], [0, 1, 2]];}{%
\[
\mathit{I1} := [[], \,[2], \,[0], \,[1], \,[0, \,1], \,[0, \,2], 
\,[1, \,2], \,[0, \,1, \,2]]
\]
}
\end{maplelatex}

\end{maplegroup}
\begin{maplegroup}
\begin{mapleinput}
\mapleinline{active}{1d}{#I123:=GenerateAdjacentIntervals(P,rel);}{%
}
\end{mapleinput}

\end{maplegroup}
\begin{maplegroup}
\begin{mapleinput}
\mapleinline{active}{1d}{CH:=[op(zip((x,y)->x=y,P,CHp)),[0,1]=[],[1,2]=[],P=1];}{%
}
\end{mapleinput}

\mapleresult
\begin{maplelatex}
\mapleinline{inert}{2d}{CH := [0 = 1, 1 = 0, 2 = 1, [0, 1] = [], [1, 2] = [], [0, 1, 2] =
1];}{%
\[
\mathit{CH} := [0=1, \,1=0, \,2=1, \,[0, \,1]=[], \,[1, \,2]=[], 
\,[0, \,1, \,2]=1]
\]
}
\end{maplelatex}

\end{maplegroup}
\begin{maplegroup}
\begin{mapleinput}
\mapleinline{active}{1d}{Con:=ConnectionMatrices(P,rel,CH,var,"full_notation");}{%
}
\end{mapleinput}

\mapleresult
\begin{maplelatex}
\mapleinline{inert}{2d}{Con := [[matrix([[0, 1, 0], [0, 0, 0], [0, 1, 0]])]];}{%
\[
\mathit{Con} :=  \left[  \!  \left[  \!  \left[ 
{\begin{array}{rrr}
0 & 1 & 0 \\
0 & 0 & 0 \\
0 & 1 & 0
\end{array}}
 \right]  \!  \right]  \!  \right] 
\]
}
\end{maplelatex}

\end{maplegroup}
\begin{maplegroup}
Now we consider $[0,1]$ as one point which we call $1$:

\end{maplegroup}
\begin{maplegroup}
\begin{mapleinput}
\mapleinline{active}{1d}{Q:=[1,2];}{%
}
\end{mapleinput}

\mapleresult
\begin{maplelatex}
\mapleinline{inert}{2d}{Q := [1, 2];}{%
\[
Q := [1, \,2]
\]
}
\end{maplelatex}

\end{maplegroup}
\begin{maplegroup}
This means, we take $C(1):=C_0(1)\oplus C_1(1) = \F_2\oplus\F_2$ and
$C(2):=C_0(2)\oplus C_1(2)=0\oplus\F_2$:

\end{maplegroup}
\begin{maplegroup}
\begin{mapleinput}
\mapleinline{active}{1d}{Cp:=[1=[1,1],2=1];}{%
}
\end{mapleinput}

\mapleresult
\begin{maplelatex}
\mapleinline{inert}{2d}{Cp := [1 = [1, 1], 2 = 1];}{%
\[
\mathit{Cp} := [1=[1, \,1], \,2=1]
\]
}
\end{maplelatex}

\end{maplegroup}
\begin{maplegroup}
\begin{mapleinput}
\mapleinline{active}{1d}{Crel:=[[2,1]];}{%
}
\end{mapleinput}

\mapleresult
\begin{maplelatex}
\mapleinline{inert}{2d}{Crel := [[2, 1]];}{%
\[
\mathit{Crel} := [[2, \,1]]
\]
}
\end{maplelatex}

\end{maplegroup}
\begin{maplegroup}
The homology {\sc Conley} index of the black box $1$ is $CH(1)=0$:

\end{maplegroup}
\begin{maplegroup}
\begin{mapleinput}
\mapleinline{active}{1d}{CCH:=[1=[],Q=1];}{%
}
\end{mapleinput}

\mapleresult
\begin{maplelatex}
\mapleinline{inert}{2d}{CCH := [1 = [], [1, 2] = 1];}{%
\[
\mathit{CCH} := [1=[], \,[1, \,2]=1]
\]
}
\end{maplelatex}

\end{maplegroup}
\begin{maplegroup}
The procedure {\tt CConnectionMatrices} reads off the missing {\sc
Conley} indices from its third argument {\tt Cp}, so the complete data
will be $CH(1)=0$, $CH(2)=CH_0(2)\oplus CH_1(2)=0\oplus \F_2$ and
$CH(\{1,2\})=CH_0(\{1,2\})\oplus CH_1(\{1,2\})=0\oplus \F_2$:

\end{maplegroup}
\begin{maplegroup}
\begin{mapleinput}
\mapleinline{active}{1d}{CCon:=CConnectionMatrices(Q,Crel,Cp,CCH,var,"full_notation");}{%
}
\end{mapleinput}

\mapleresult
\begin{maplelatex}
\mapleinline{inert}{2d}{CCon := [[matrix([[1, 0], [0, 0]])], [matrix([[1, 0], [1, 0]])]];}{%
\[
\mathit{CCon} := [[ \left[ 
{\begin{array}{rr}
1 & 0 \\
0 & 0
\end{array}}
 \right] ], \,[ \left[ 
{\begin{array}{rr}
1 & 0 \\
1 & 0
\end{array}}
 \right] ]]
\]
}
\end{maplelatex}

\end{maplegroup}
\begin{maplegroup}
Since $C(0)$ was chosen to be different from $CH(0)$ the corresponding
diagonal entry of the $C$-connection matrix does not vanish:

\end{maplegroup}
\begin{maplegroup}
\begin{mapleinput}
\mapleinline{active}{1d}{map(ConnectionSubmatrix,CCon,Q,[1,1],map(rhs,Cp),var);}{%
}
\end{mapleinput}

\mapleresult
\begin{maplelatex}
\mapleinline{inert}{2d}{[[matrix([[1]])], [matrix([[1]])]];}{%
\[
[[ \left[ 
{\begin{array}{r}
1
\end{array}}
 \right] ], \,[ \left[ 
{\begin{array}{r}
1
\end{array}}
 \right] ]]
\]
}
\end{maplelatex}

\end{maplegroup}
\begin{maplegroup}
\begin{mapleinput}
\mapleinline{active}{1d}{map(ConnectionSubmatrix,CCon,Q,[2,2],map(rhs,Cp),var);}{%
}
\end{mapleinput}

\mapleresult
\begin{maplelatex}
\mapleinline{inert}{2d}{[[matrix([[0]])], [matrix([[0]])]];}{%
\[
[[ \left[ 
{\begin{array}{r}
0
\end{array}}
 \right] ], \,[ \left[ 
{\begin{array}{r}
0
\end{array}}
 \right] ]]
\]
}
\end{maplelatex}

\end{maplegroup}
\begin{maplegroup}
\begin{mapleinput}
\mapleinline{active}{1d}{map(ConnectionSubmatrix,CCon,Q,[1,2],map(rhs,Cp),var);}{%
}
\end{mapleinput}

\mapleresult
\begin{maplelatex}
\mapleinline{inert}{2d}{[[matrix([[0]])], [matrix([[0]])]];}{%
\[
[[ \left[ 
{\begin{array}{r}
0
\end{array}}
 \right] ], \,[ \left[ 
{\begin{array}{r}
0
\end{array}}
 \right] ]]
\]
}
\end{maplelatex}

\end{maplegroup}
\begin{maplegroup}
For the specified data there exist two compatible $C$-connection
matrices. The second connection matrix would imply a connecting orbit
from $2$ to $1$:

\end{maplegroup}
\begin{maplegroup}
\begin{mapleinput}
\mapleinline{active}{1d}{map(ConnectionSubmatrix,CCon,Q,[2,1],map(rhs,Cp),var);}{%
}
\end{mapleinput}

\mapleresult
\begin{maplelatex}
\mapleinline{inert}{2d}{[[matrix([[0]])], [matrix([[1]])]];}{%
\[
[[ \left[ 
{\begin{array}{r}
0
\end{array}}
 \right] ], \,[ \left[ 
{\begin{array}{r}
1
\end{array}}
 \right] ]]
\]
}
\end{maplelatex}

\end{maplegroup}
\begin{maplegroup}
In particular, computations with $C$-connection matrices allow one to
consider {\sc Morse} decompositions with {\sc Morse} sets having
trivial homology {\sc Conley} indices, i.e.\  with $CH(p)=0$.

\end{maplegroup}
%
\POP

\section{Conclusion}
We noticed during the implementation that the braid condition usually appearing in the definition of $(C)$-connection matrices is, due to {\sc Franzosa}'s result in Prop.~\ref{complexbraid}, automatically satisfied and does not impose any further restrictions on the matrices. This led us to the simplified Definition \ref{conn}. We adapted the classical definitions to the modern language of derived categories and introduced triangles and octahedra. Since sufficiently many of the examples studied in dynamical systems are symmetric, we gave a precise definition of a symmetric $C$-connection matrix and studied an example with {\em non}-hyperbolic equilibria, where the symmetry group acts non-trivially on the homology {\sc Conley} indices of the non-hyperbolic equilibria.

The purpose of the package {\tt conley} \cite{conley} is to automatize long computations and free the user from lots of technical details, by this allowing him to deal with large (generic) examples, which can hardly be processed by hand. We hope that the package eliminates some of the difficulties one might encounter while entering this, in our opinion, very interesting area of dynamical systems.

Larger and more interesting examples than in Section \ref{beispiele} will be studied in future joint work with {\sc Maier-Paape}. There {\tt conley}'s ability to compute transition matrices will be another essential ingredient.

%
\POP

\PUSH{conley.bbl}%
\providecommand{\bysame}{\leavevmode\hbox to3em{\hrulefill}\thinspace}
\providecommand{\MR}{\relax\ifhmode\unskip\space\fi MR }

\providecommand{\href}[2]{#2}

\POP


\begin{thebibliography}{MPMW}

\bibitem[BR]{BR}
Mohamed Barakat and Daniel Robertz, \emph{{$\mathtt{homalg}$ -- {A}
  meta-package for homological algebra}}, submitted. {\tt
  arXiv:math.AC/0701146} and (\href{http://wwwb.math.rwth-aachen.de/homalg}{\tt
  http://wwwb.math.rwth-aachen.de/homalg}).

\bibitem[BR07]{conley}
\bysame, \emph{$\mathtt{conley}$ project}, 2006-2007,
  (\href{http://wwwb.math.rwth-aachen.de/conley}{\tt
  http://wwwb.math.rwth-aachen.de/conley}).

\bibitem[Con78]{Con}
Charles Conley, \emph{Isolated invariant sets and the {M}orse index}, CBMS
  Regional Conference Series in Mathematics, vol.~38, American Mathematical
  Society, Providence, R.I., 1978. \MR{MR511133 (80c:58009)}

\bibitem[Fra86]{Fra86}
Robert~D. Franzosa, \emph{Index filtrations and the homology index braid for
  partially ordered {M}orse decompositions}, Trans. Amer. Math. Soc.
  \textbf{298} (1986), no.~1, 193--213. \MR{MR857439 (88a:58121)}

\bibitem[Fra88]{Fra88}
\bysame, \emph{The continuation theory for {M}orse decompositions and
  connection matrices}, Trans. Amer. Math. Soc. \textbf{310} (1988), no.~2,
  781--803. \MR{MR973177 (90g:58111)}

\bibitem[Fra89]{Fra89}
\bysame, \emph{The connection matrix theory for {M}orse decompositions}, Trans.
  Amer. Math. Soc. \textbf{311} (1989), no.~2, 561--592. \MR{MR978368
  (90a:58149)}

\bibitem[GM03]{GM}
Sergei~I. Gelfand and Yuri~I. Manin, \emph{Methods of homological algebra},
  second ed., Springer Monographs in Mathematics, Springer-Verlag, Berlin,
  2003. \MR{MR1950475 (2003m:18001)}

\bibitem[Mis95]{Mis95}
Konstantin Mischaikow, \emph{Conley index theory}, Dynamical systems
  (Montecatini Terme, 1994), Lecture Notes in Math., vol. 1609, Springer,
  Berlin, 1995, pp.~119--207. \MR{MR1374109 (97a:58109)}

\bibitem[MM02]{Mis}
Konstantin Mischaikow and Marian Mrozek, \emph{Conley index theory}, Handbook
  of Dynamical Systems III: Towards Applications (B.~Fiedler, G.~Iooss, and
  N.~Kopell, eds.), North-Holland, 2002.

\bibitem[MPMW]{MMW}
S.~Maier-Paape, K.~Mischaikow, and T.~Wanner, \emph{{Structure of the attractor
  of the Cahn-Hilliard equation on the square}}, International J. Bifurcation
  Chaos, 70 pp., to appear.

\end{thebibliography}


\edoc
